\documentclass[10pt,journal,]{IEEEtran}
\usepackage{times}
\usepackage{hyperref}
\usepackage{url}
\usepackage{amsmath}
\usepackage{amssymb}
\usepackage{graphicx}
\usepackage{color}
\usepackage[tight,footnotesize]{subfigure}

\newcommand{\norm}[1]{\left\Vert#1\right\Vert}
\newcommand{\abs}[1]{\left\vert#1\right\vert}
\newcommand\vect[1]{{\bf#1}}
\newcommand\matr[1]{{\bf#1}}
\newcommand{\Real}{\mathbb R}
\newcommand\RR[1]{\mathbb{R}^{#1}}

\newcommand{\Kappa}{\mathcal K}
\newcommand{\Proj}{\mathcal P}
\newcommand{\shrink}{\mathcal S}

\newtheorem{thm}{Theorem}[section]

\newtheorem{lem}[thm]{Lemma}

\newtheorem{defn}[thm]{Definition}

\usepackage{algorithmic}
\usepackage{algorithm}

\DeclareMathOperator{\argmin}{argmin}
\DeclareMathOperator{\sgn}{sgn}

\ifCLASSOPTIONcompsoc
  \usepackage[nocompress]{cite}
\else
  \usepackage{cite}
\fi

\title{Tradeoffs between Convergence Speed and Reconstruction Accuracy in Inverse Problems}

\usepackage{authblk}

\author[1]{Raja Giryes}
\author[2]{Yonina C. Eldar}
\author[3]{Alex M. Bronstein}
\author[4]{Guillermo Sapiro}

\affil[1]{School of Electrical Engineering\\
  Tel Aviv University University\\
  Tel Aviv, Israel 69978 
\url{raja@tauex.tau.ac.il}}
\affil[2]{  Electrical Engineering Department\\
  Technion - IIT\\
  Haifa, Israel, 32000\\
\url{yonina@ee.technion.ac.il}}
\affil[3]{  Computer Science Department\\
  Technion - IIT\\
  Haifa, Israel, 32000\\
\url{bron@cs.technion.ac.il}}
\affil[4]{    Electrical and Computer Engineering Department\\
  Duke University\\
  Durham, NC, 27708\\
\url{guillermo.sapiro@duke.edu}}

\begin{document}

\maketitle

\begin{abstract}
Solving inverse problems with iterative algorithms is popular, especially for large data. Due to time constraints, the number of possible iterations is usually limited, potentially affecting the achievable accuracy. Given an error one is willing to tolerate, an important question is whether it is possible to modify the original iterations to obtain faster convergence to a minimizer achieving the allowed error without increasing the computational cost of each iteration considerably.
Relying on recent recovery techniques developed for settings in which the desired signal belongs to some low-dimensional set, we show that using a coarse estimate of this set may lead to faster convergence at the cost of an additional reconstruction error related to the accuracy of the set approximation.
Our theory ties to recent advances in sparse recovery, compressed sensing, and deep learning. Particularly, it may provide a possible explanation to the successful approximation of the $\ell_1$-minimization solution by neural networks with  layers representing iterations, as practiced in the learned iterative shrinkage-thresholding algorithm (LISTA).
\end{abstract}

\section{Introduction}

Consider the setting in which we want to recover a vector $\vect{x} \in \RR{d}$ from  linear measurements
\begin{eqnarray}
\label{eq:meas}
\vect{y} = \matr{M}\vect{x} + \vect{e},
\end{eqnarray}
where $\matr{M} \in \RR{m \times d}$ is the measurement matrix and $\vect{e} \in \RR{d}$ is additive noise. 
This setup appears in many fields including statistics (e.g., regression), image processing (e.g., deblurring and super-resolution), and medical imaging (e.g., CT and MRI), to name just a few.

Often the recovery of $\vect{x}$ from $\vect{y}$ is an ill-posed problem. For example, when $\matr{M}$ has fewer rows than columns ($m < n$), rendering \eqref{eq:meas} an underdetermined linear system of equations. In this case, it is impossible to recover $\vect{x}$ without introducing additional assumptions on its structure. A popular strategy is to assume that $\vect{x}$ resides in a low dimensional set $\Kappa$, e.g., sparse vectors \cite{Bruckstein09From, Candes05Decoding, Elad10Sparse,  Eldar15Sampling} or a Gaussian Mixture Model (GMM) \cite{Yu12Solving}. The natural by-product minimization problem then becomes
\begin{eqnarray}
\label{eq:min_const}
\min_{\vect{x}} \norm{\vect{y} - \matr{M}\vect{x}}_2^2  & \text{s.t.} & \vect{x} \in \Kappa.
\end{eqnarray} 
This can be reformulated in an unconstrained form as
\begin{eqnarray}
\label{eq:min_unconst}
\min_{\vect{x}} \norm{\vect{y} - \matr{M}\vect{x}}_2^2 + \lambda f(\vect{x}),
\end{eqnarray} 
where $\lambda$ is a regularization parameter and $f(\cdot)$ is a cost function related to the set $\Kappa$. For example, if $\Kappa = \left\{ \vect{x} \in \RR{d} : \norm{\vect{x}}_0 \le k\right\}$ is the set of $k$-sparse vectors, then a natural choice is $f(\cdot) = \norm{\cdot}_0$ or its convex relaxation $f(\cdot) = \norm{\cdot}_1$.

A popular technique for solving $\eqref{eq:min_const}$ and $\eqref{eq:min_unconst}$ is using iterative programs such as proximal methods \cite{Combettes11Proximal, Bottou16Optimization} that include the iterative shrinkage-thresholding algorithm (ISTA) \cite{Beck09Fast, Elad07Coordinate, Daubechies04iterative} and the alternating direction method of multipliers (ADMM) \cite{Boyd11Distributed, Gabay76dual}. This strategy is particularly useful for large dimensions $d$. 

Many applications impose  time constraints, which limit the number of computations that can be performed to recover $\vect{x}$ from the measurements. One way to minimize time and computations is to reduce the number of iterations without increasing the computational cost of each iteration. A different approach is to use momentum methods \cite{Wilson16Lyapunov} or random projections \cite{Pilanci15Iterative, Pilanci16Iterative, Pilanci15Newton, Tu17Breaking} to accelerate convergence.
Another alternative is to keep the number of iterations fixed while reducing the cost of each iteration.
For example, since the complexity of iterative methods rely, among other things, on $m$, a common technique to save computations is to sub-sample the measurements $\vect{y}$, removing ``redundant information,'' to an amount that still allows reconstruction of $\vect{x}$. A series of recent works \cite{Bottou08Tradeoffs, Shalev08SVM, Daniely13More, Bruer15Designing, Chandrasekaran13Computational, Oymak15Sharp} suggest that by obtaining more measurements one can benefit from simple efficient methods that cannot be applied with a smaller number of measurements. 

In \cite{Bottou08Tradeoffs} the generalization properties of large-scale learning systems have been studied showing a tradeoff between the number of measurements and the target approximation. The work in \cite{Shalev08SVM} showed how it is possible to make the run-time of SVM optimization decrease as the size of the training data increases. In \cite{Daniely13More}, it is shown that the problem of supervised learning of halfspaces over 3-sparse vectors with trinary values $\{-1, 1, 0\}$ may be solved with efficient algorithms only if the number of training examples exceeds a certain limit. Similar phenomena are encountered  in the context of sparse recovery, where efficient algorithms are guaranteed to reconstruct the sparsest vector only if the number of samples is larger than a certain quantity \cite{Eldar12Compressed, Elad10Sparse, Foucart13Mathematical}. 
In \cite{Chandrasekaran13Computational} it was shown that by having a larger number of training examples it is possible to design more efficient optimization problems by projecting onto simpler sets. This idea is further studied in \cite{Bruer15Designing} by changing the amount of smoothing applied in convex optimization. In \cite{Oymak15Sharp} the authors show that more measurements may allow increasing the step-size in the projected gradient algorithm (PGD) and thus accelerating its convergence.

While these works studied a tradeoff between convergence speed and the number of available measurements, this paper takes a different route. 
Consider the case in which due to time constraints we need to stop the iterations before we achieve the desired reconstruction accuracy. For the original algorithm, this can result in the recovery being very far from the optimum.
An important question is whether we can modify the original iterations (e.g., those dictated by the shrinkage or ADMM techniques), such that the method convergences to an improved solution with fewer iterations without adding complexity to them.
This introduces a tradeoff between the recovery error we are willing to tolerate and the computational cost. As we demonstrate, this goes beyond the trivial relationship between the approximation error and the number of iterations that exists for various iterative methods \cite{Beck09Fast}. 

Such a tradeoff is experimentally demonstrated by the success of {\it learned} ISTA (LISTA) \cite{Gregor10Learning} for sparse recovery with $f(\cdot) = \norm{\cdot}_1$. This technique learns a neural network with only several layers, where each layer is a modified version of the ISTA iteration.\footnote{ISTA and its variants is one of the most powerful optimization techniques for sparse coding.} It achieves virtually the same accuracy as the original ISTA using one to two orders of magnitude less iterations. The acceleration of iterative algorithms with neural networks is not unique only to the sparse recovery problem and $f(\cdot) = \norm{\cdot}_1$. This behavior was demonstrated for other models such as the analysis cosparse and low-rank matrix models \cite{Sprechmann15Learning}, Poisson noise \cite{Remez15Picture}, acceleration of Eulerian fluid simulation \cite{Tompson16Accelerating}, and feature learning \cite{Andrychowicz16Learning}.
However, a proper theoretical justification to this phenomena is still lacking.

{\bf Contribution.} In this work, we provide theoretical foundations elucidating the tradeoff between the allowed minimization error and the number of simple iterations used  for solving inverse problems. We formally show that if we allow a certain reconstruction error in the solution, then it is possible to change iterative methods by modifying the linear operations applied in them such that each iteration has the same complexity as before but the number of steps required to attain a certain error is reduced. 

Such a tradeoff seems natural when working with real data, where both the data and the assumed models are noisy or approximate; searching for the {\it exact} solution of an optimization problem, where all the variables are affected by measurement or model noise  may be an unnecessary use of valuable computational resources.
We formally prove this relation for iterative projection algorithms.
Interestingly, a related tradeoff exists also in the context of sampling theory, where by allowing some error in the reconstruction we may use fewer samples and/or quantization levels \cite{Kipnis16Sampling}.
We argue that the tradeoff we analyze may explain the smaller number of iterations required in LISTA compared to ISTA.

%An important observation is that there is a difference between the {\em reconstruction error}, i.e., the error relative to the true signal, and the minimization error in the {\em objective value}, i.e., the error in the objective of the minimization problem. We may be violating significantly the latter while doing very well in the former, which is what we would like to reduce.

Parallel efforts to our work also provide justification for the success of  LISTA. In \cite{Moreau16Adaptive}, the fast convergence of LISTA is justified by connecting between the convergence speed and the factorization of the Gram matrix of $\matr{M}$. In \cite{Xin16Maximal}, the convergence speed of ISTA and LISTA is analyzed using the restricted isometry property (RIP) \cite{Candes05Decoding}, showing that LISTA may reduce the RIP, which leads to faster convergence. A relation between LISTA and approximate message passing (AMP) strategies is drawn in \cite{Borgerding16Onsager}.

Our paper differs from previous contributions in three main points: (i) it goes beyond the case of standard LISTA  with sparse signals and considers variants that apply to general low-dimensional models; (ii) our theory relies on the concept of inexact projections and their relation to the tradeoff between convergence-speed and recovery accuracy, which differs significantly from other attempts to explain the success of LISTA; and (iii) besides exploring LISTA, we provide acceleration strategies to other programs such as model-based compressed sensing \cite{Eldar12Compressed} and sparse recovery with side-information. 

{\bf Organization.} This paper is organized as follows. In Section~\ref{sec:prelim} we present preliminary notation and definitions, and describe the ISTA, LISTA and PGD techniques. Section~\ref{sec:PGD_new_theory} introduces a new theory for PGD for non-convex cones. 
Section~\ref{sec:iter_proj_tradeoofs} shows how it is possible to tradeoff between convergence speed and reconstruction accuracy by introducing the inexact projected gradient descent (IPGD) method using spectral compressed sensing \cite{Duarte13Spectral} as a motivating example. The reconstruction error of IPGD is analyzed as a function of the iterations in Section~\ref{sec:IPGD_theory}. 
Section~\ref{sec:examples} discusses the relation between our theory and model-based compressed sensing \cite{Baraniuk10Model} and sparse recovery with side information \cite{Friedlander12Recovering, Khajehnejad09Weighted, Vaswani10Modified, Weizman15Compressed}. Section~\ref{sec:LIPGD} proposes a LISTA version of IPGD, the learned IPGD (LIPGD), and demonstrates its usage in the task of image super-resolution.
Section~\ref{sec:cs_DL} relates the approximation of minimization problems studied here with neural networks  and deep learning, providing a theoretical foundation for the success of LISTA and suggesting a ``mixture-model'' extension of this technique. Section~\ref{sec:conc} concludes the paper.

\section{Preliminaries and Background}
\label{sec:prelim}

Throughout this paper, we use the following notation. We write 
$\norm{\cdot}$ for the Euclidian norm for vectors and the spectral norm for matrices, $\norm{\cdot}_1$ for the $\ell_1$
    norm that sums the absolute values of a vector and $\norm{\cdot}_0$
    for the $\ell_0$ pseudo-norm, which counts the number of non-zero elements in a vector. The conjugate transpose of $\matr{M}$ is denoted by 
$\matr{M}^*$ and the orthogonal projection onto the set $\Kappa$ by $\Proj_{\Kappa}$. The original unknown vector is denoted by $\vect{x}$, the given measurements by $\vect{y}$, the measurement matrix by $\matr{M}$ and the system noise by $\vect{e}$. The $i$th entry of a vector $\vect{v}$ is denoted by $\vect{v}[i]$. The sign function $\sgn(\cdot)$ equals $1$, $-1$ or $0$ if its input is positive, negative or zero respectively. The $d$-dimensional $\ell_2$-ball of radius $r$ is denoted by $\mathbb{B}^d_r$. For balls of radius $1$, we omit the subscript and just write $\mathbb{B}^d$.

\subsection{Iterative shrinkage-thresholding algorithm (ISTA)}

%% TODO: add refernce from Miki's book about condition for ISTA step-size
A popular iterative technique for minimizing \eqref{eq:min_unconst} is ISTA.  Each of its iterations is composed of a gradient step with step size $\mu$, obeying $\frac{1}{\mu} \ge \norm{\matr{M}}$ to ensure convergence \cite{Beck09Fast}, followed by a proximal mapping $\shrink_{f,\lambda}(\cdot)$ of the function $f$, defined as
\begin{eqnarray}
\shrink_{f,\lambda }(\vect{v}) = \argmin_{\vect{z}} \frac{1}{2}\norm{\vect{z} - \vect{v}} + \lambda f(\vect{z}),
\end{eqnarray}
where $\lambda$ is a parameter of the mapping. 
The resulting ISTA iteration can be written as
\begin{eqnarray}
\label{eq:ISTA}
\vect{z}_{t+1} = \shrink_{f,\mu\lambda}\left(\vect{z}_t + \mu \matr{M}^*(\vect{y} - \matr{M}\vect{z}_t)\right),
\end{eqnarray}
where $\vect{z}_t$ is an estimate of $\vect{x}$ at iteration $t$. Note that the step size $\mu$ multiplies the parameter of the proximal mapping.

The proximal mapping has a simple form for many functions $f$. For example, when $f(\cdot) = \norm{\cdot}_1$, it is an element-wise shrinkage function,
\begin{eqnarray}
\label{eq:soft_thresh}
\shrink_{\ell_1,\lambda}(\vect{v})[i] = \sgn\left(\vect{v}[i]\right)\max(0, \abs{\vect{v}[i]} - \lambda).
\end{eqnarray}
Therefore, the advantage of ISTA is that its iterations require only the application of matrix multiplications and then a simple non-linear function.
Nonetheless, the main drawback of ISTA is the large number of iterations that is typically required for convergence. 

\begin{figure}[t]
\begin{center}
{\includegraphics[width=0.45\textwidth]{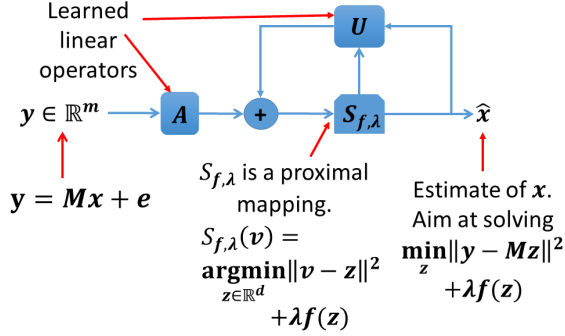}
}
\end{center}
\vspace{-0.15in}
\caption{The LISTA scheme.}
\label{fig:LISTA_scheme}
\end{figure}

Many acceleration techniques have been proposed to speed up convergence of ISTA (see \cite{Beck09Fast, Beck09FastGradient, Elad07Coordinate, Nesterov83Method, Treister12Multilevel, Goldstein09Split, Su14Differential, Machart12Optimal, Aujol15Stability, Wibisono16variational} as a partial list of such works).
A prominent strategy is LISTA, which has the same structure as ISTA but with different linear operations in \eqref{eq:ISTA}. Empirically, it is observed that it is able to attain a solution very close to that of ISTA with a significantly smaller fixed number of iterations $T$. The LISTA iterations are given by\footnote{We present the more general version \cite{Sprechmann15Learning} that can be used for any signal model and not only for sparsity.}
\begin{eqnarray}
\label{eq:LISTA}
\vect{z}_{t+1} = \shrink_{f,\lambda}\left(\matr{A}\vect{y} + \matr{U}\vect{z}_t\right),
\end{eqnarray}
where $\matr{A}$, $\matr{U}$ and $\lambda$ are learned from a set of training examples by back-propagation with the objective being the $\ell_2$-distance between the final ISTA solution and the LISTA one (after $T$ iterations) \cite{Gregor10Learning}. 
Other minimization objectives may be used, e.g., training LISTA to minimize \eqref{eq:min_unconst} directly \cite{Sprechmann15Learning}. Notice that LISTA has a structure of a recurrent neural network as can be seen in Fig.~\ref{fig:LISTA_scheme}.

While other acceleration techniques for ISTA have been proposed together with a thorough theoretical analysis, the powerful LISTA method has been introduced without mathematical justification for its success. 
In this work, we focus on the PGD algorithm, whose iterations are almost identical to the ones of ISTA but with an orthogonal projection instead of a proximal mapping. We propose an acceleration technique for it, which is very similar to the one of LISTA, accompanied by a theoretical analysis. 

\subsection{Projected gradient descent (PGD)}

The PGD iteration is given by
\begin{eqnarray}
\label{eq:PGD}
\vect{z}_{t+1} = \Proj_{\Kappa}\left(\vect{z}_t + \mu \matr{M}^*(\vect{y} - \matr{M}\vect{z}_t)\right),
\end{eqnarray}
where $\Proj_{\Kappa}$ is an orthogonal projection onto a given set $\Kappa$.
%Note that these iterations are almost identical to the ones of ISTA but with just one difference: instead of a proximal mapping, an orthogonal projection is performed onto a given set $\Kappa$.
For example, if $\Kappa$ is the $\ell_1$-ball then $\Proj_{\Kappa}$ is simply soft thresholding with a  value that varies depending on the projected vector \cite{Duchi08Efficient}. Note the similarity to the proximal mapping in ISTA with $f$ as the $\ell_1$-norm, which is also the soft thresholding operation but with a fixed threshold \eqref{eq:soft_thresh}. This similarity is not unique to the $\ell_1$-norm case but happens also for other types of $f$ such as the $\ell_0$ pseudo-norm and the nuclear norm.
%and $f$ is the $\ell_1$-norm, then in PGD the projection $\Proj_\Kappa$  is simply soft threshodling with a  value that varies depending on the projected vector \cite{Duchi08Efficient} and in LISTA the proximal projection is  also soft-thresholding but with a fixed threshold .
The step size $\mu$ is assumed to be constant for the sake of simplicity, as in \eqref{eq:ISTA}. In both methods it may vary between  iterations.

PGD is a generalization of the iterative hard thresholding (IHT) algorithm, which was developed for $\Kappa$ being the set of sparse vectors \cite{Blumensath09Iterative}. 
This important method has been analyzed in various works. For example, for standard sparsity in  \cite{Blumensath09Iterative}, for sparsity patterns that belong to a certain model in  \cite{Baraniuk10Model}, for a general union of subspaces in \cite{Blumensath11Sampling}, for nonlinear measurements in \cite{Beck13Sparsity}, 
and more recently in \cite{Oymak15Sharp} for a set of the form $$\Kappa = \left\{ \vect{z} \in \RR{d} : f(\vect{z}) \le R  \right\}.$$ The formulation \eqref{eq:PGD} generalizes the special cases above. For example, if $f(\cdot) = \norm{\cdot}_0$ and $R$ is the sparsity level then we have the IHT method from \cite{Blumensath09Iterative}; when $f$ counts the number of non-zeros of only certain sparsity patterns, which are bounded by $R$, we have the model-based IHT of  \cite{Baraniuk10Model}. PGD may also be applied to non-linear inverse problems \cite{Yang16Sparse, Oymak15Sharp}.

Theorem~\ref{thm:PGD_error_t} below provides convergence guarantees on PGD (it is the noiseless version of Theorem 1.2 in \cite{Oymak15Sharp}).
Before presenting the result, we introduce several properties of the set $\Kappa$ and some basic  lemmas.
\begin{defn}[Descent set and tangent cone]
The descent set of the function $f$ at a point $\vect{x}$ is defined as
\begin{eqnarray}
D_f(\vect{x}) = \left\{ \vect{h} \in \RR{d} : f(\vect{x} + \vect{h}) \le f(\vect{x})  \right\}.
\end{eqnarray}
The tangent cone $C_f(\vect{x})$ at a point $\vect{x}$ is the conic hull of $D_f(\vect{x})$, i.e., the smallest closed cone $C_f(\vect{x})$ satisfying  $D_f(\vect{x}) \subseteq C_f(\vect{x})$.
\end{defn}
For concise writing, below we denote $D_f(\vect{x})$ and $C_f(\vect{x})$ as $D$ and $C$, respectively.
\begin{lem}[Lemma 6.2 in \cite{Oymak15Sharp}]
\label{lem:P_c_sup}
Let $\vect{v}\in \Real^d$ and $C \subset \Real^d$ be a closed cone. Then 
\begin{eqnarray}
\norm{\Proj_C(\vect{v})} = \sup_{\vect{u} \in C \cap \mathbb{B}^d} \vect{u}^* \vect{v}.
\end{eqnarray}
\end{lem}

\begin{lem}
\label{lem:Kappa_D_rel}
If for $\vect{x} \in \RR{d}$, $\Kappa = \{ \vect{z} \in \RR{d} : f(\vect{z}) \le f(\vect{x}) \} \subset \RR{d}$ is a closed set, then for all $\vect{v} \in \RR{d}$, 
\begin{eqnarray}
\Proj_\Kappa(\vect{x} + \vect{v} ) -\vect{x} = \Proj_{\Kappa - \left\{\vect{x} \right\}}(\vect{v}) = \Proj_{D} (\vect{v} ).
\end{eqnarray}
\end{lem}
{\it Proof:} From the definition of the descent cone we have
 $D = \left\{\vect{h} \in \RR{d} :  f(\vect{h} + \vect{x}) \le f(\vect{x}) \right\} =   \left\{\vect{z} - \vect{x} :  f(\vect{z}) \le f(\vect{x}) \right\} = \left\{\vect{z} - \vect{x} : \vect{z} \in \Kappa \right\} = \Kappa - \left\{\vect{x} \right\}$,
where  the second equality follows from a simple change of variables, and the last ones from the definitions of the set $\Kappa$ and the Minkowski difference. Therefore, projecting onto $D$ is equivalent to a projection onto $\Kappa - \left\{\vect{x} \right\}$.
 \hfill $\Box$ 

\begin{lem}[Lemma~6.4 in \cite{Oymak15Sharp}]
\label{lem:C_D_rel}
Let ${D}$ and $C$ be a nonempty and closed set and cone, respectively, such that $0 \in D$ and $D \subseteq C$. Then for all $\vect{v} \in \RR{d}$ 
\begin{equation}
\norm{\Proj_D(\vect{v})} \le \kappa_f \norm{\Proj_C(\vect{v})},
\end{equation} 
where $\kappa_f = 1$ if $D$ is a convex set and $\kappa_f = 2$ otherwise.
\end{lem}

We now introduce the convergence rate provided in \cite{Oymak15Sharp} for PGD. For brevity, we present only its noiseless version.

\begin{thm}[Noiseless version of Theorem~1.2 in \cite{Oymak15Sharp}]
\label{thm:PGD_error_t}
Let $\vect{x} \in \RR{d}$, $f: \RR{d} \rightarrow \Real$ be a proper function, $\Kappa = \left\{ \vect{z} \in \RR{d} : f(\vect{z}) \le f(\vect{x}) \right\}$, 
$C = C_f(\vect{x})$ the tangent cone of the function $f$ at point $\vect{x}$, $\matr{M} \in \RR{m \times d}$ and $\vect{y} = \matr{M}\vect{x}$ a vector containing $m$ linear measurements. Assume we are using PGD with $\Kappa$  to recover $\vect{x}$ from $\vect{y}$. Then the estimate $\vect{z}_t$ at the $t$th iteration  (initialized with $\vect{z}_0 = 0$) obeys
\begin{eqnarray}
\norm{\vect{z}_t - \vect{x}} \le (\kappa_f \rho(C))^t\norm{\vect{x}},
\end{eqnarray}
where $\kappa_f$ is defined in Lemma~\ref{lem:C_D_rel}, and
\begin{eqnarray}
\label{eq:rho_C}
\rho(C) = \rho(\mu, \matr{M}, f, \vect{x}) = \sup_{\vect{u}, \vect{v} \in C \cap \mathbb{B}^d} \vect{u}^*\left(\matr{I} - \mu\matr{M}^*\matr{M} \right)\vect{v},
\end{eqnarray}
is the convergence rate of PGD.
\end{thm}

\subsection{Gaussian mean width}
\label{sec:mean_width_def}

When $\matr{M}$ is a random matrix with i.i.d. Gaussian distributed entries $\mathcal{N}(0,1)$, it has been shown in \cite{Oymak15Sharp} that the convergence rate 
$\rho(C)$ is tightly related to the dimensionality of the set (model) $\vect{x}$ resides in. A very useful expression for measuring the ``intrinsic dimensionality'' of sets is the (Gaussian) mean width. 
\begin{defn}[Gaussian mean width] The Gaussian mean width of a set $\Upsilon$ is defined as
\begin{eqnarray}
\label{eq:omega}
\omega(\Upsilon) = E[\sup_{\vect{v} \in \Upsilon \cap \mathbb{B}^d} \langle \vect{g}, \vect{v} \rangle ], & \vect{g} \sim \mathcal{N}(0,\matr{I}). 
\end{eqnarray}
\end{defn}
Two variants of this measure are generally used. The {\em cone Gaussian mean width},
$\omega_C = \omega(C)$,
which measures the dimensionality of the tangent cone $C = C_f(\vect{x})$; and the {\em set Gaussian mean width}, $\omega_\Kappa = \omega(\Kappa - \Kappa)$, which is related directly to the set $\Kappa$ through its Minkowski difference $\Kappa - \Kappa = \left\{ \vect{z} - \vect{v} : \vect{z},\vect{v} \in \Kappa  \right\}$.
%\begin{defn}[Cone Gaussian mean width] The Gaussian mean width of the tangent cone $C=C_f(\vect{x})$ is defined as
%\end{defn}
%\begin{defn}[Set Gaussian mean width] The Gaussian mean width of a set $\Kappa$ is defined as \end{defn}} (intersected with the unit $\ell_2$-ball $\mathbb{B}^d$)
The cone Gaussian mean width relies on both the set $\Kappa$ (through $f$) and a specific target point $\vect{x}$, while the set Gaussian mean width considers only $\Kappa$. On the other hand, the dependence of $\omega_C$ on $\Kappa$ is indirect via the descent set at the point $\vect{x}$. %, where the latter treats $\Kappa$ directly via its Minkowski difference. 
There is a series of works, which developed convergence and reconstruction guarantees for various methods based on $\omega_C$ \cite{Chandrasekaran12Convex, Amelunxen14Living, Oymak15Sharp}, and others that rely on $\omega_\Kappa$ \cite{Plan13Robust, Tirer17Generalizing}. The first ($\omega_C$) is mainly employed in the case of convex functions $f$, which are used to  relax the non-convex set in which $\vect{x}$ resides. In this setting, often $D$ is convex and $\vect{x} \in \Kappa$.

As an example of $\omega_C$ consider the case in which $\Kappa$ is the $\ell_1$-ball and $\vect{x}$ is a $k$-sparse vector. Then $\omega_C \simeq \sqrt{2 k \log(d/k)}$.
If we add constraints on $\vect{x}$ such as having a tree structure, i.e., belonging to the set
\begin{eqnarray}
\label{eq:Kappa_hat_tree}
\hat{\Kappa} = \{  \vect{z}\in \RR{d} : \norm{\vect{z}}_0 \le \norm{\vect{x}}_0  \text{\&} ~ \vect{z} ~ \text{obeys a tree structure} \},
\end{eqnarray}
where an entry may be non-zero only if its parent node is non-zero,
then the value of $\omega_C$ does not change. Although the definition of $\omega_\Kappa$ is very similar to $\omega_C$, it yields different results. For $\Kappa$ the set of $k$-sparse vectors $\omega_\Kappa = O(\sqrt{k\log(d/k)})$, while for \eqref{eq:Kappa_hat_tree}, $\omega_{\hat{\Kappa}} = O(\sqrt{k})$. The first result is similar to the expression of $\omega_C$ for the $\ell_1$ ball with a $k$-sparse vector $\vect{x}$, yet, the second provides a better measure of the set $\hat{\Kappa}$ in \eqref{eq:Kappa_hat_tree}.

\subsection{PGD convergence rate and the cone Gaussian mean width}
\label{sec:conv_Gaussian_width}

In \cite{Oymak15Sharp}, it has been shown that the smaller $\omega_C$, the faster the convergence. More specifically, if $m$ is very close to $\omega_C$, then we may apply PGD with a step-size $\mu =\frac{1}{(\sqrt{d}+\sqrt{m})^2}  \simeq \frac{1}{d}$  and have  a convergence rate of (Theorem~2.4 in \cite{Oymak15Sharp})
\begin{eqnarray}
\label{eq:rho_omegaC1}
\rho(C) = 1- O\left(\frac{\sqrt{m} - \omega_C }{m+d}\right).
\end{eqnarray} 
If $\omega_C$ is smaller than $\sqrt{m}$ by a certain constant factor, then we may apply PGD with a larger step size $\mu \simeq \frac{1}{m}$, which leads to improved convergence (Theorem~2.2 in \cite{Oymak15Sharp})
\begin{eqnarray}
\label{eq:rho_omegaC2}
\rho(\Kappa) = O\left(\frac{\omega_C}{\sqrt{m}}\right).
\end{eqnarray}

 These relationships rely on the fact that with larger $m$ the eigenvalues of $\matr{I} - \mu \matr{M}^T\matr{M}$ (after projection onto $C$, see \eqref{eq:rho_C}) are better positioned such that it is possible to improve convergence by increasing $\mu$. 

Both \eqref{eq:rho_omegaC1} and \eqref{eq:rho_omegaC2}  set a limit on the minimal value $m$ for which PGD iterations  converge to $\vect{x}$, namely $m=O(\omega_C^2)$. This implies that $m \gtrapprox 2 k \log(d/k)$ (bigger than approximately $2 k \log(d/k)$)
for $\Kappa$ as the $\ell_1$-ball and a $k$-sparse vector $\vect{x}$. This is known to be a tight condition. See more examples for this relationship between $\omega_C$ and $m$ in \cite{Chandrasekaran12Convex}.

%\rg{In Section~\ref{sec:PGD_new_theory} we derive theory for PGD, which is related to .}

The connections in \eqref{eq:rho_omegaC1} and \eqref{eq:rho_omegaC2} between $\rho(C)$ and $\omega_C$ are not unique only to the case that $\matr{M}$ is a random Gaussian matrix. Similar relationships hold for many other types of matrices \cite{Oymak15Sharp}.

\section{PGD Theory based on the Projection Set}
\label{sec:PGD_new_theory}

While Theorem~\ref{thm:PGD_error_t} covers many sets $\Kappa$, there are interesting examples that are not included in it such as the set of $k$-sparse vectors corresponding to $f$ being the $\ell_0$ pseudo-norm, which is not a proper function. Even if we ignore this condition and try to use the result of Theorem~\ref{thm:PGD_error_t} in the case that $\matr{M}$ is a random Gaussian matrix we face a problem. Using the relationship between $\rho(C)$ and the Gaussian mean width $\omega_C$ in \eqref{eq:rho_omegaC1} and \eqref{eq:rho_omegaC2}, and the fact that in this case $\omega_C= \sqrt{d}$, 
we get the condition $m >d$. This demand on $m$ is inferior to existing theory that in this scenario guarantees convergence with $m = O(k \log(d/k)$ \cite{Blumensath09Iterative}.

One way to overcome this problem is by considering the convex-hull of the set of $k$-sparse vectors with bounded $\ell_2$ norm. In this case  $\omega_C = O(\sqrt{k\log(d/k)})$ \cite{Plan13Robust}. However, as mentioned above, guarantees for PGD exist for the $k$-sparse case without a bound on the $\ell_2$-norm.

A similar phenomenon also occurs with the set of sparse vectors with a tree structure $\hat{\Kappa}$ (see \eqref{eq:Kappa_hat_tree}),
where again $\omega_C = \sqrt{d}$ implying $m = O(d)$. Yet, from the work in \cite{Baraniuk10Model}, we know that in this setting it is sufficient to choose $m = O(k)$. Note that for $\hat{\Kappa}$, the set Gaussian mean width is $\omega_{\hat{\Kappa}} = \sqrt{k}$. If we would have relied on it instead of on $\omega_C$ in the bound for the required size of $m$,  it would have coincided with \cite{Baraniuk10Model}.

In order to address these deficiencies in the convergence rate, we provide a variant of Theorem~\ref{thm:PGD_error_t} that relies on the set $\Kappa$ directly through the Minkowski difference $\Kappa - \Kappa$ in lieu of $C$. 
 For simplicity we present only the noiseless case but the extension to the noisy setting can be performed using the strategy in \cite{Oymak15Sharp}. 
\begin{thm}
\label{thm:PGD_error_K}
Let $\vect{x} \in \Kappa$, 
$\Kappa \subset \RR{d}$ be a closed cone,
$\matr{M} \in \RR{m \times d}$ and $\vect{y} = \matr{M}\vect{x}$ a vector containing $m$ linear measurements. Assume we are using PGD with $\Kappa$  to recover $\vect{x}$ from $\vect{y}$. Then the estimate $\vect{z}_t$ at the $t$th iteration  (initialized with $\vect{z}_0 = 0$) obeys
\begin{eqnarray}
\norm{\vect{z}_t - \vect{x}} \le (\kappa_\Kappa \rho(\Kappa))^t\norm{\vect{x}},
\end{eqnarray}
where $\kappa_\Kappa=1$ if $\Kappa$ is convex and $\kappa_\Kappa = 2$ otherwise, and
\begin{eqnarray}
\rho(\Kappa) = \rho(\mu, \matr{M}, \Kappa) = \hspace{-0.1in} \sup_{\vect{u}, \vect{v} \in (\Kappa - \Kappa) \cap \mathbb{B}^d} \hspace{-0.1in} \vect{u}^*\left(\matr{I} - \mu\matr{M}^*\matr{M} \right)\vect{v},
\end{eqnarray}
is the convergence rate of PGD.
\end{thm}

{\it Proof:}
We repeat similar steps to the ones in the proof of Theorem~1.2 in \cite{Oymak15Sharp}. 

We start by noting that the PGD error at iteration $t+1$ is,
\begin{eqnarray}
\label{eq:PGD_error_def}
\hspace{-0.27in}  \norm{\vect{z}_{t+1} - \vect{x}} &=& \norm{\Proj_{\Kappa}\left(\vect{z}_t + \mu \matr{M}^*(\vect{y} - \matr{M}\vect{z}_t)\right) - \vect{x}}
 \\ \nonumber  \hspace{-0.22in} &=&  \norm{\Proj_{D}\left(\left( \matr{I} - \mu \matr{M}^*\matr{M}\right)(\vect{z}_t - \vect{x}) \right) },
\end{eqnarray}
where the last inequality is due to Lemma~\ref{lem:Kappa_D_rel} and the fact that $\vect{y} = \matr{M}\vect{x}$. Since $\Kappa$ is a closed cone, also the Minkowski difference $\Kappa - \Kappa$ is a closed cone. Moreover, $D \subset \Kappa-\Kappa$ as $\vect{x} \in \Kappa$. Thus, following Lemma~\ref{lem:C_D_rel} we have 
\begin{eqnarray}
\label{eq:PGD_step_error_proof_1}
&& \hspace{-0.27in} \norm{\vect{z}_{t+1} - \vect{x}}  
\le  \kappa_\Kappa\norm{\Proj_{\Kappa-\Kappa}\left(\left( \matr{I} - \mu \matr{M}^*\matr{M}\right)(\vect{z}_t - \vect{x}) \right) }
\\ \nonumber && \hspace{-0.27in}
\le \hspace{-0.1in} \sup_{\substack{\vect{v} \in \Real^d \text{ s.t. } \\   \norm{\Proj_{\Kappa}(\vect{v})- \vect{x}} 
\le \norm{\vect{z}_t - \vect{x}}}}  \hspace{-0.1in} 
\kappa_\Kappa\norm{\Proj_{\Kappa-\Kappa}\left(\left( \matr{I} - \mu \matr{M}^*\matr{M}\right)(\Proj_\Kappa(\vect{v}) - \vect{x}) \right) }
\\ \nonumber && \hspace{-0.27in}
\le \hspace{-0.1in} \sup_{\substack{\vect{v} \in \Real^d \text{ s.t. } \\   \norm{\Proj_D(\vect{v}- \vect{x})} 
\le \norm{\vect{z}_t - \vect{x}}}}  \hspace{-0.1in} 
\kappa_\Kappa\norm{\Proj_{\Kappa-\Kappa}\left(\left( \matr{I} - \mu \matr{M}^*\matr{M}\right)(\Proj_D(\vect{v} - \vect{x})) \right) },
\end{eqnarray}
where the second inequality is due to the fact that $\vect{z}_t$ is of the form $\Proj_\Kappa(\vect{v})$ for some vector $\vect{v}$, and the last inequality follows from Lemma~\ref{lem:Kappa_D_rel}. 

Noticing that the constraint $\norm{\Proj_D(\vect{v}- \vect{x})} 
\le \norm{\vect{z}_t - \vect{x}}$ is equivalent to $\vect{v}- \vect{x} \in D \cap  \mathbb{B}^d_{\norm{\vect{z}_t - \vect{x}}}$ (where $ \mathbb{B}^d_{\norm{\vect{z}_t - \vect{x}}}$ is the $\ell_2$-ball of radius $\norm{\vect{z}_t - \vect{x}}$) and 
using the relation $D \subset \Kappa-\Kappa$ leads to 
\begin{eqnarray}
\label{eq:PGD_step_error_proof_2}
&& \hspace{-0.27in} \norm{\vect{z}_{t+1} - \vect{x}}  
\\ \nonumber && \hspace{-0.27in}
\le  \sup_{\vect{v} \in (\Kappa-\Kappa) \cap \mathbb{B}^d}
\kappa_\Kappa\norm{\Proj_{\Kappa-\Kappa}\left(\left( \matr{I} - \mu \matr{M}^*\matr{M}\right)\vect{v} \right) }\norm{\vect{z}_t - \vect{x}}
\\ \nonumber && \hspace{-0.27in}
\le  \sup_{\vect{v}, \vect{u} \in (\Kappa-\Kappa) \cap \mathbb{B}^d}
\kappa_\Kappa\norm{\vect{u}^*\left(\left( \matr{I} - \mu \matr{M}^*\matr{M}\right)\vect{v} \right) }\norm{\vect{z}_t - \vect{x}},
\end{eqnarray}
where the last inequality follows from Lemma~\ref{lem:P_c_sup}.
Using the definition of $\rho(\Kappa)$ and applying the inequality in \eqref{eq:PGD_step_error_proof_2} recursively leads to the desired result.  
 \hfill $\Box$ 

When $\matr{M}$ is a random Gaussian matrix, the relationships in \eqref{eq:rho_omegaC1} and \eqref{eq:rho_omegaC2} hold with $\rho(\Kappa)$ and $\omega_{\Kappa}$ replacing $\rho(C)$ and $\omega_C$ respectively. This implies that we need $m = O(\omega_\Kappa^2)$ for convergence. This result is in line with the conditions on $m$ that appear in previous works for $k$-sparse vectors \cite{Blumensath09Iterative}, for which $\omega_\Kappa = O(\sqrt{k\log(d/k)})$,  and for sparse vectors with tree structure \cite{Baraniuk10Model}, where $\omega_\Kappa = O(\sqrt{k})$.

As discussed in Section~\ref{sec:mean_width_def}, the measure $\omega_\Kappa$ is related directly to the set $\Kappa$ (may be non-convex) in which $\vect{x}$ resides. 
Thus, it provides a better measure for the complexity of $\Kappa$ when it is unbounded or has some specific structure as is the case for sparsity with tree structure \cite{Tirer17Generalizing}. In such settings, Theorem~\ref{thm:PGD_error_K} should be favored over Theorem~\ref{thm:PGD_error_t}.

Notice that if $\vect{x} \in \Kappa$, then we have $D \subset \Kappa - \Kappa$. Thus, in the settings that $D$ is convex and $\vect{x} \in \Kappa$, we have $C = D \subset \Kappa - \Kappa$ implying that $\omega_C \le \omega_\Kappa$; when $\matr{M}$ is random Gaussian, this also implies $\rho(C) \le \rho(\Kappa)$. Therefore, in this scenario Theorem~\ref{thm:PGD_error_t} has an advantage over Theorem~\ref{thm:PGD_error_K}.

\section{Inexact projected gradient descent (IPGD)}
\label{sec:iter_proj_tradeoofs}

It may happen that the function $f$ or the set $\Kappa$ are too loose for describing $\vect{x}$. Instead, we may select a set $\hat{\Kappa}$ that better characterizes $\vect{x}$ and therefore leads to a smaller $\omega$, resulting in faster convergence. This improvement can be very significant; smaller $\omega$ both improves the convergence rate and allows using  a larger step-size (see Section~\ref{sec:conv_Gaussian_width}). 

For example, consider the case of a $k$-sparse vector $\vect{x}$, whose sparsity pattern obeys a tree structure. 
%(an entry may be non-zero only if its parent node is non-zero). 
If we ignore the structure in $\vect{x}$ and choose $f$ as the $\ell_0$ or $\ell_1$ norms, then the mean widths are $\omega_{\Kappa} = O(k \log(d/k))$ \cite{Plan13Robust} and $\omega_{C} \simeq 2(k \log(d/k))$ \cite{Chandrasekaran12Convex} respectively. However, if we take this structure into account and use the set of $k$-sparse vectors with tree structure (see \eqref{eq:Kappa_hat_tree}),
then $\omega_{\hat{\Kappa}} = O(k)$ \cite{Tirer17Generalizing}.  As mentioned above, this improvement may be significant especially when $m$ is very close to $\omega_{\Kappa}$.  
Such an approach was taken in the context of model-based compressed sensing \cite{Baraniuk10Model}, where it is shown that faster convergence is achieved by projecting onto the set of $k$-sparse vectors with tree structure instead of the standard $k$-sparse set. 

A related study \cite{Yu12Solving} showed that it is enough to use a small number of Gaussians to represent all the patches in natural images instead of using a dictionary that spans a much larger union of subspaces. This work relied on Gaussian Mixture Models (GMM), whose mean width scales proportionally to the number of Gaussians used, which is significantly smaller than the mean width of the sparse model.

\subsection{Inexact projection}
A difficulty often encountered is that the projection onto $\hat\Kappa$, which may even be unknown, is more complex to implement than the projection onto $\Kappa$. The latter can be easier to project onto but provides a lower convergence rate. 

Thus, in this work we introduce a technique that compromises between the reconstruction error and convergence speed by using  
PGD with an inexact ``projection'' that projects onto a set that is approximately as small as $\hat{\Kappa}$ but yet is as computationally efficient as the projection onto $\Kappa$. In this way, the computational complexity of each projected gradient descent iteration remains the same while the convergence rate becomes closer to that of the more complex PGD with a projection onto $\hat{\Kappa}$.

The ``projection'' we propose is composed of a simple operator $p$ (e.g., a linear  or an element-wise function) and the projection onto $\Kappa$, $\Proj_{\Kappa}$, such that it introduces only a slight distortion into $\vect{x}$. In particular, we require the following:

\subsubsection{The projection condition for convex \texorpdfstring{$\Kappa$}{K}}
If $\Kappa$ is convex, then we require 
\begin{eqnarray}
\label{eq:Kp_eps_ineq}
\norm{\vect{x} - \Proj_{\Kappa}\left( p(\vect{x})\right)}  \le \epsilon\norm{\vect{x}}.
\end{eqnarray} 
Due to Lemma~\ref{lem:Kappa_D_rel}, this is equivalent to
\begin{eqnarray}
\label{eq:Dp_eps_ineq}
\norm{\Proj_{D}\left(\vect{x} - p(\vect{x})\right)} \le \epsilon\norm{\vect{x}}.
\end{eqnarray} 
From the fact that $\norm{\Proj_{D}\left(\vect{x} - p(\vect{x})\right)} \le \norm{\vect{x} - p(\vect{x})}$, it is sufficient that
\begin{eqnarray}
\label{eq:p_eps_ineq}
\norm{\vect{x} - p(\vect{x})} \le \epsilon\norm{\vect{x}},
\end{eqnarray}
to ensure \eqref{eq:Kp_eps_ineq}. 
Examples for projections that satisfy condition \eqref{eq:Kp_eps_ineq} are given hereafter in sections~\ref{sec:IPGD_desc} and \ref{sec:sparse_side}.

\subsubsection{The projection condition for non-convex \texorpdfstring{$\Kappa$}{K}}
In the case that $\Kappa$ is non-convex, we require
\begin{equation}
\label{eq:P_K_pv_x_px_epsilon}
\norm{\Proj_\Kappa\left( p \vect{v} - p\vect{x} \right) - \Proj_\Kappa\left( p \vect{v} - \vect{x} \right)}  \le \epsilon\norm{\vect{x}}, \forall \vect{v} \in \Real^d.
\end{equation}
Due to Lemma~\ref{lem:Kappa_D_rel} and a simple change of variables, \eqref{eq:P_K_pv_x_px_epsilon} is equivalent to 
\begin{equation}
\label{eq:P_D_pv_x_px_epsilon}
\norm{\Proj_D\left( p \vect{v} - \vect{x} + p\vect{x} \right) - \Proj_D\left( p \vect{v} \right)}  \le \epsilon\norm{\vect{x}}, \forall \vect{v} \in \Real^d,
\end{equation}
which by another simple change of variables is the same as
\begin{equation}
\label{eq:P_D_pv_x_pv_px_epsilon}
\norm{\Proj_D\left( p \vect{v} - \vect{x} \right) - \Proj_D\left( p \vect{v} - p\vect{x}\right)}  \le \epsilon\norm{\vect{x}}, \forall \vect{v} \in \Real^d.
\end{equation} 
An example for a projection that satisfies condition \eqref{eq:P_K_pv_x_px_epsilon} is provided in Section~\ref{sec:tree}.

\begin{figure*}[t]
\begin{center}
{
{
\subfigure[x2 redundant DCT dictionary with sparsity $k=2$.]{\includegraphics[width=0.48\textwidth]{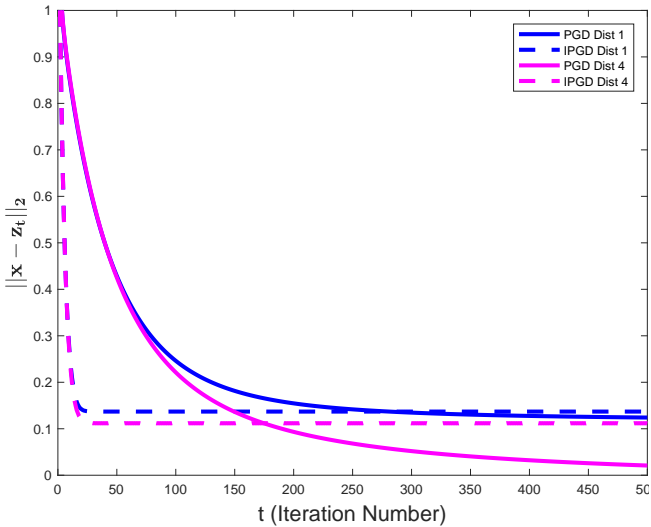}
\label{fig:spect_error_2redundant}}
\hspace{0.05in}
\subfigure[x4 redundant DCT dictionary with sparsity $k=4$.]{\includegraphics[width=0.48\textwidth]{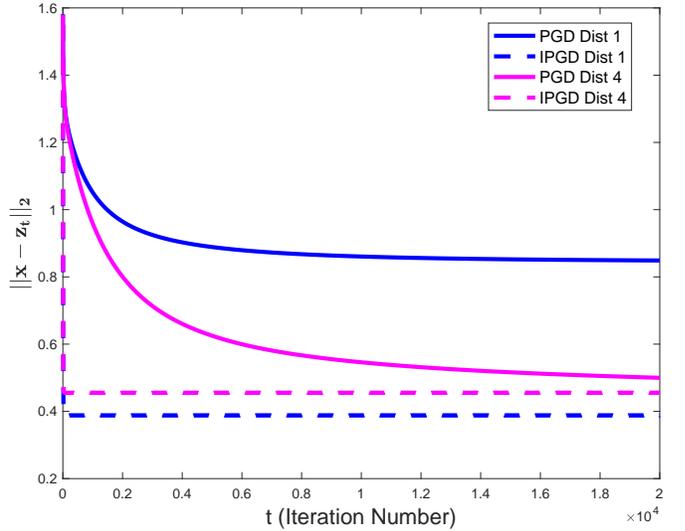}
\label{fig:spect_error_4redundant}}
}
}
\end{center}
\caption{Reconstruction error as a function of iterations for sparse recovery with a dictionary with high coherence between neighboring atoms. The sparse representation in the dictionary is generated such that there are three correlated neighboring atoms close to each other with location distance 1 or 4.
PGD is applied with $\Kappa$ being the $\ell_1$ ball. IPGD is used with the same $\Kappa$ and $p$ being a non-linear function that for a given vector keeps at most only one dominant entry in every neighborhood of fixed size (zeroing the smaller values). 
This shows that IPGD may accelerate convergence compared to PGD and in some cases (right figure) even achieve lower recovery error.
}
\label{fig:spect_error}
\end{figure*}

\subsection{Inexact PGD}
\label{sec:IPGD_desc} 

Plugging the inexact projection into the PGD step results in the proposed {\it inexact} PGD (IPGD) iteration (compare to \eqref{eq:PGD})
\begin{eqnarray}
\label{eq:IPGD}
\vect{z}_{t+1} = \Proj_{\Kappa}\left(p\left(\vect{z}_t\right) + \mu p\left( \matr{M}^*(\vect{y} - \matr{M}\vect{z}_t)\right)\right).
\end{eqnarray}
To motivate this algorithm consider the problem of spectral compressed sensing \cite{Duarte13Spectral}, in which one wants to recover a sparse representation in a dictionary that has high local coherence. It has been shown that if the non-zeros in the representation are far from each other then it is easier to obtain good recovery \cite{Candes14Towards}. 

Let $\matr{M}$ be a two times redundant DCT dictionary and $\tilde{\vect{x}}$ be a $k$-sparse vector, with sparsity $k=2$, of dimension $d =128$, such that the minimal distance (with respect to the location in the vector) between non-zero neighboring coefficients in it is greater than $5$ (indices) and the value in each non-zero coefficient is generated from the normal distribution. We construct the vector $\vect{x}$ by adding to $\tilde{\vect{x}}$ random Gaussian values with zero mean and variance $\sigma^2 = 0.05$ at the neighboring coefficients of each non-zero entry in $\tilde{\vect{x}}$ with (location) distance $1$ or $4$ (two different experiments).

As mentioned above, a better reconstruction is achieved by estimating $\tilde{\vect{x}}$ from $\matr{M}\tilde{\vect{x}}$ than by estimating $\vect{x}$ from $\matr{M}\vect{x}$ due to the highly correlated columns in $\matr{M}$. A common practice to improve the recovery in such a case is to force the recovery algorithm to select a solution with separated coefficients. In our context it is simply using the IPGD with a projection onto the $\ell_1$ ball and $p(\cdot)$ that keeps at most only one dominant entry (in absolute value) in every neighborhood of size $5$ in a given representation by zeroing out the other values. 
The operator $p$ causes an error in the model (with $\epsilon\simeq 0.05\sqrt{2} \simeq 0.1$) and therefore reaches a slightly higher final error than PGD with projection onto the $\ell_1$ ball.  Compared to PGD, IPGD projects onto a simpler set with a smaller Gaussian mean width, thus, attaining faster convergence at the first iterations, where the approximation error is still significantly larger than $\epsilon$ as can be seen in Fig.~\ref{fig:spect_error_2redundant}. When the coherence is larger (as in the case of added coefficients at distance $1$), the advantage of IPGD over PGD is more significant.
 
%\rg{such as in the case of added entries at distance $1$, where there is a large coherence between the atoms in the representation (compared to the case of distance $4$ that has much lower coherence). For another such example,}

In some cases IPGD may even attain a lower final recovery error compared to PGD. For example, consider the case of $\matr{M}$ being a four times redundant DCT dictionary and $\vect{x}$ generated as above but with $k=4$. %be equal to $\tilde{\vect{x}}$ with additional random Gaussian values with zero mean and variance $\sigma^2 = 0.05$ at the neighboring coefficients at distance $1$ and $2$ of each non-zero entry in $\tilde{\vect{x}}$.  
Due to the larger redundancy in the dictionary, the coherence is larger in this case. Thus, the recovery of $\vect{x}$ is harder.
Here, PGD with a projection onto the $\ell_1$ ball converges slower and reaches a large error due to the high correlations between the atoms. Using IPGD with the $\ell_1$ ball and a projection $p(\cdot)$ that keeps at most only one dominant entry (in absolute value) in every neighborhood of size $5$, leads both to faster convergence and better final accuracy as can be seen in Fig.~\ref{fig:spect_error_4redundant}. %ADD an illustration of how the signal looks like
 
\section{IPGD Convergence Analysis} 
 \label{sec:IPGD_theory}

We turn to analyze the performance of IPGD.
For simplicity of the discussion, we analyze the convergence of this technique only for a linear operator $p$ and the noiseless setting, i.e., $\vect{e} = 0$.  The extension to other types of operators and the noisy case is straightforward by arguments similar to those used in \cite{Oymak15Sharp} for treating the noise term and other classes of matrices. 

We present two theorems on the convergence of IPGD. The first result provides a bound in terms of $\rho(C)$ (i.e., depends on $\omega_C$ if $\matr{M}$ is a random Gaussian matrix) for the case that $D$ is convex corresponding to $\kappa_f=1$ in Theorem~\ref{thm:PGD_error_t};  the second provides a bound in terms of $\rho(K)$ (i.e., depends on $\omega_\Kappa$ if $\matr{M}$ is a random Gaussian matrix) when $\Kappa$ is a closed cone but not necessarily convex. 
The proofs of both theorems are deferred to appendices~\ref{sec:IPGD_error_t_proof} and \ref{sec:IPGD_error_K_proof}.

\begin{thm}
\label{thm:IPGD_error_t}
Let $\vect{x} \in \RR{d}$, $f: \RR{d} \rightarrow \Real$ be a proper function, $\Kappa = \left\{ \vect{z} \in \RR{d} : f(\vect{z}) \le f(\vect{x}) \right\}$, 
$D = D_{f}(\vect{x})$ and $C = C_{f}(\vect{x})$ the descent set and the tangent cone of the function $f$ at point $\vect{x}$ respectively, $p(\cdot)$ a linear operator satisfying \eqref{eq:Dp_eps_ineq},
$\matr{M} \in \RR{m \times d}$ and $\vect{y} = \matr{M}\vect{x}$ a vector containing $m$ linear measurements. Assume we are using IPGD with $\Kappa$ and $p$  to recover $\vect{x}$ from $\vect{y}$ and that $D$ is convex. Then the estimate $\vect{z}_t$ at the $t$th iteration  (initialized with $\vect{z}_0 = 0$) obeys
\begin{eqnarray*}
 \norm{\vect{z}_t - \vect{x}} 
 \le \left(\left(\rho_p(C)\right)^t +  \frac{1-\left( \rho_p({C})  \right)^t}{1- \rho_p({C})   }(2+\rho_p(C))\epsilon  \right)\norm{\vect{x}},
\end{eqnarray*}
where 
\begin{eqnarray*}
\rho_p(C) &=& \rho(\mu, \matr{M}, f, p, \vect{x}) \\ \nonumber 
&=& \hspace{-0.15in} \sup_{\vect{u},\vect{v} \in C\cap \mathbb{B}^d}  
p(\vect{u})^*\left( \matr{I} - \mu \matr{M}^*\matr{M}\right)p(\vect{v})
\end{eqnarray*}
is the ``effective convergence rate'' of IPGD for small $\epsilon$.
\end{thm}

%%%%%%%%%%%%%%%%%%%%%%%%%%%%%%%%%%%%%%%%%
%%         Non-convex case
%%%%%%%%%%%%%%%%%%%%%%%%%%%%%%%%%%%%%%%%%

\begin{thm}
\label{thm:IPGD_error_K}
Let $\vect{x} \in \Kappa$, 
$\Kappa \subset \RR{d}$ be a closed cone, $p(\cdot)$ a linear operator satisfying \eqref{eq:P_K_pv_x_px_epsilon},
$\matr{M} \in \RR{m \times d}$ and $\vect{y} = \matr{M}\vect{x}$ a vector containing $m$ linear measurements. Assume we are using IPGD with $\Kappa$ and $p$  to recover $\vect{x}$ from $\vect{y}$. Then the estimate $\vect{z}_t$ at the $t$th iteration  (initialized with $\vect{z}_0 = 0$) obeys
\begin{eqnarray}
&& \hspace{-0.5in} \norm{\vect{z}_t - \vect{x}} 
 \le \left(\left(\kappa_\Kappa\rho_p(\Kappa)\right)^t +  \frac{1-\left( \kappa_\Kappa\rho_p({\Kappa})  \right)^t}{1- \kappa_\Kappa\rho_p({\Kappa})   }\gamma  \right)\norm{\vect{x}},
\end{eqnarray}
where $\kappa_\Kappa$ and $\rho(K)$ are defined in Theorem~\ref{thm:PGD_error_K},
\begin{eqnarray}
\gamma \triangleq (2\rho(\Kappa)\kappa_\Kappa+\rho_p(\Kappa)\kappa_\Kappa+1)\epsilon,
\end{eqnarray}
and
\begin{eqnarray}
\label{eq:rho_p_K}
\rho_p(\Kappa) &=& \rho(\mu, \matr{M}, \Kappa, p) \\ \nonumber &=& \hspace{-0.15in} \sup_{\vect{u}, \vect{v} \in (\Kappa - \Kappa) \cap \mathbb{B}^d}  p(\vect{u})^*\left(\matr{I} - \mu\matr{M}^*\matr{M} \right)p(\vect{v})
\end{eqnarray}
is the ``effective convergence rate'' of IPGD for small $\epsilon$.
\end{thm}

Theorems~\ref{thm:IPGD_error_t} and \ref{thm:IPGD_error_K}  imply that if $\epsilon$ is small enough (compared to $\rho_p^t$, where $t$ is the iteration number and $\rho_p$ is defined in \eqref{eq:rho_p_K}) then IPGD has an effective convergence rate of $\rho_p=\rho_p(C)$ when $D$ is convex, and $\rho_p=\kappa_\Kappa\rho_p(\Kappa)$ in the case that $\Kappa$ is a closed cone but not necessarily convex.  Note that if $p = \matr{I}$ then $\epsilon =0$ and our results coincide with theorems~\ref{thm:PGD_error_t} and \ref{thm:PGD_error_K}.

As we shall see hereafter, for some operators $p$ the rate $\rho_p$ may be significantly smaller than $\rho(C)$ and $\rho(\Kappa)$. The smaller the set that $p$ maps to, the smaller $\rho_p$ becomes. At the same time, when $p$ maps to smaller sets it usually provides a ``coarser estimate'' and thus the approximation error $\epsilon$ in \eqref{eq:Dp_eps_ineq} and \eqref{eq:P_K_pv_x_px_epsilon} increases. 
Thus, IPGD allows us to tradeoff approximation error $\epsilon$  and  improved convergence $\rho_p$.

\begin{figure*}[t]
\begin{center}
{
{
\subfigure[Recovery error as a function of iteration number]{\includegraphics[width=0.48\textwidth]{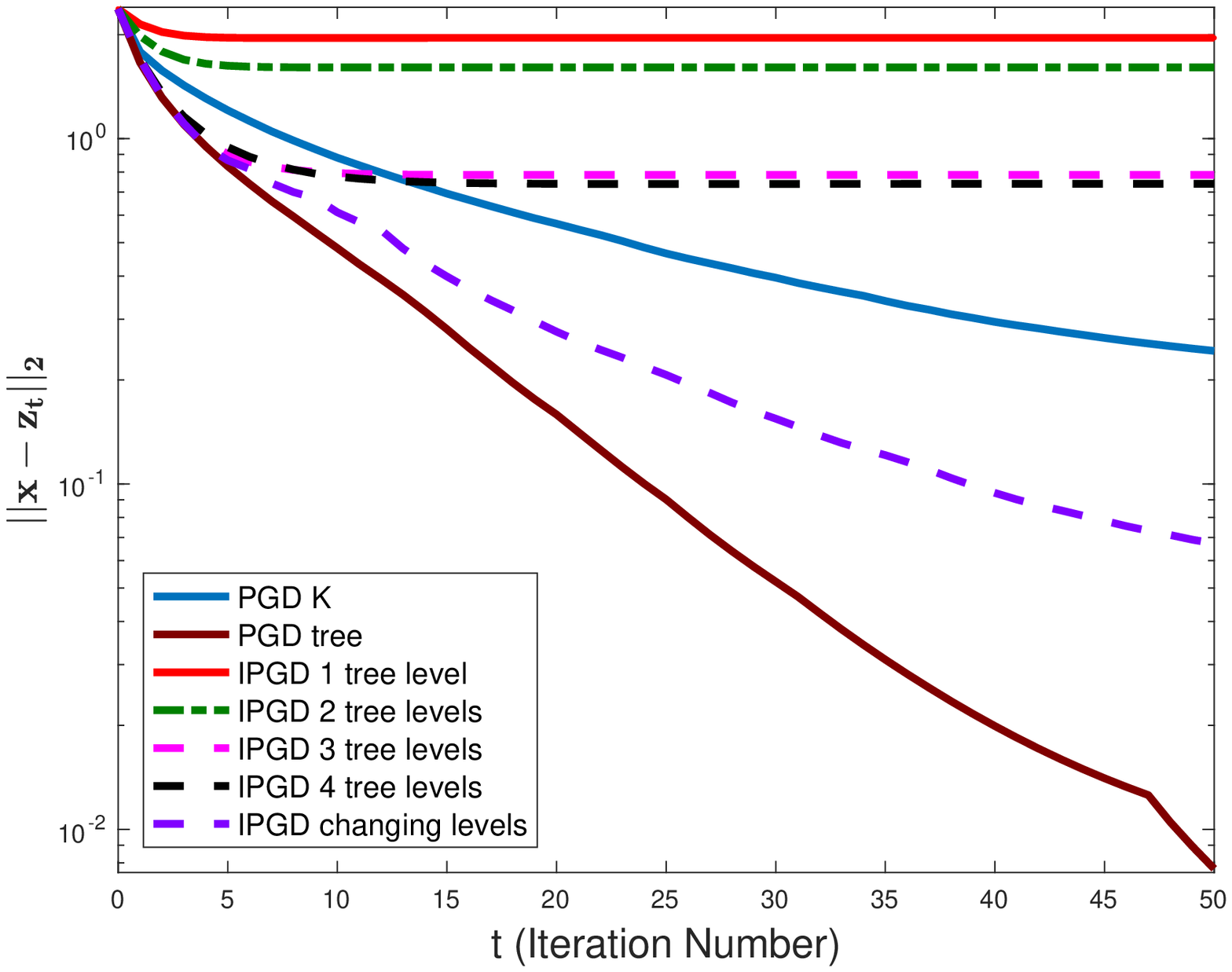}
\label{fig:tree_error}}
\hspace{0.05in}
\subfigure[Recovery error as a function of running time (in sec)]{\includegraphics[width=0.48\textwidth]{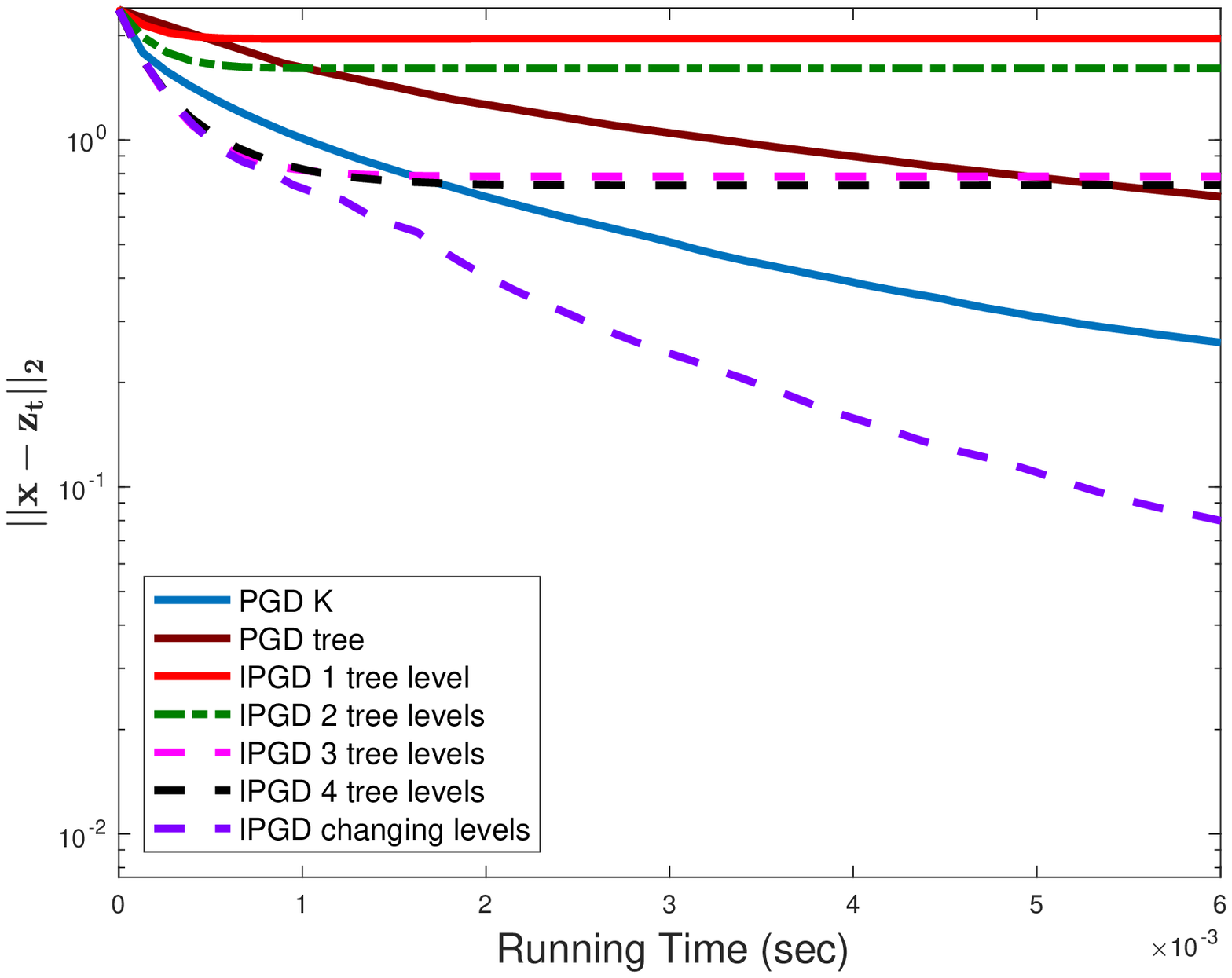}
\label{fig:tree_error_run_time}}
\subfigure[Recovery error as a function of iteration number zoomed]{\includegraphics[width=0.48\textwidth]{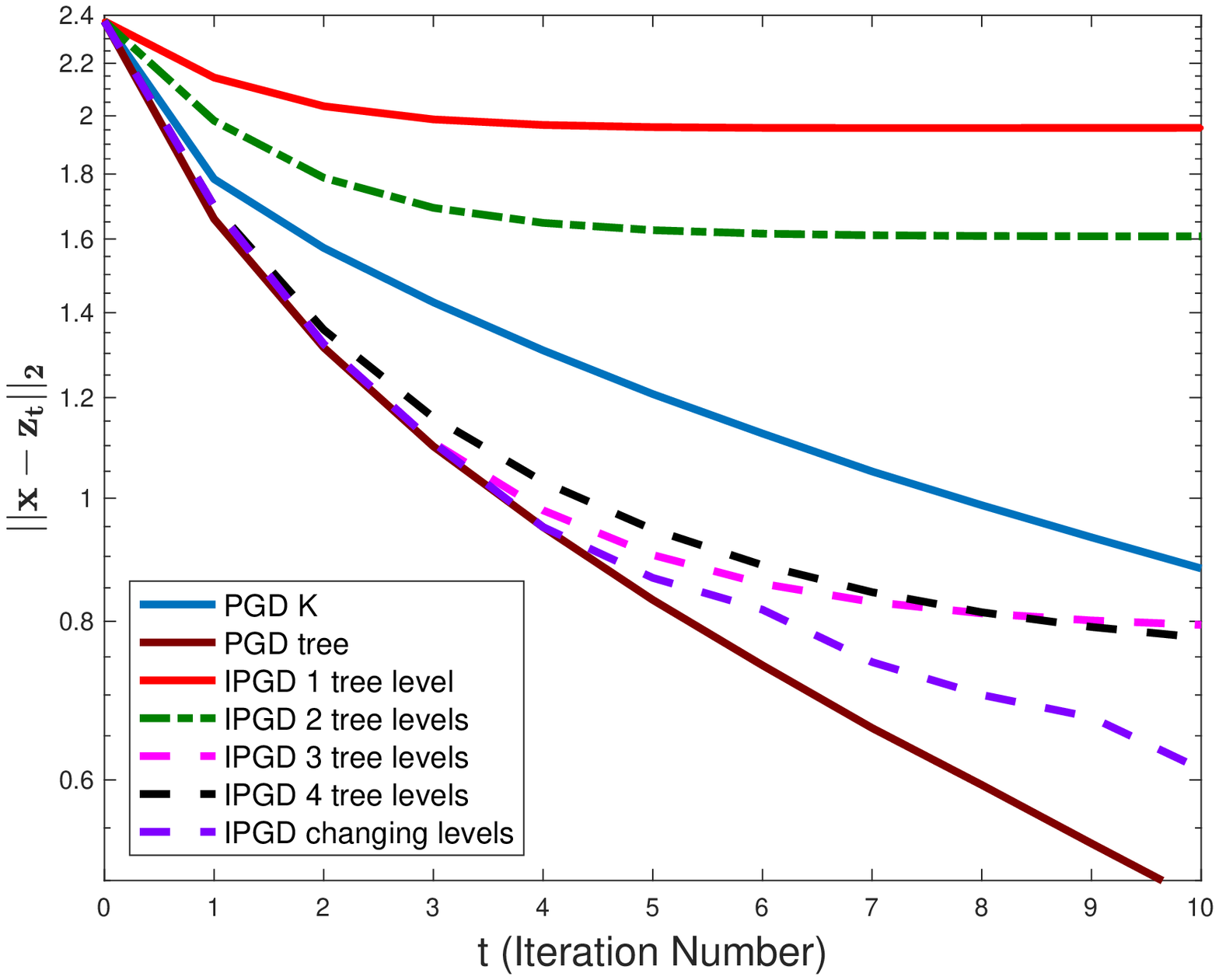}
\label{fig:tree_error_zoom}}
\hspace{0.05in}
\subfigure[Recovery error as a function of running time (in sec) zoomed]{\includegraphics[width=0.48\textwidth]{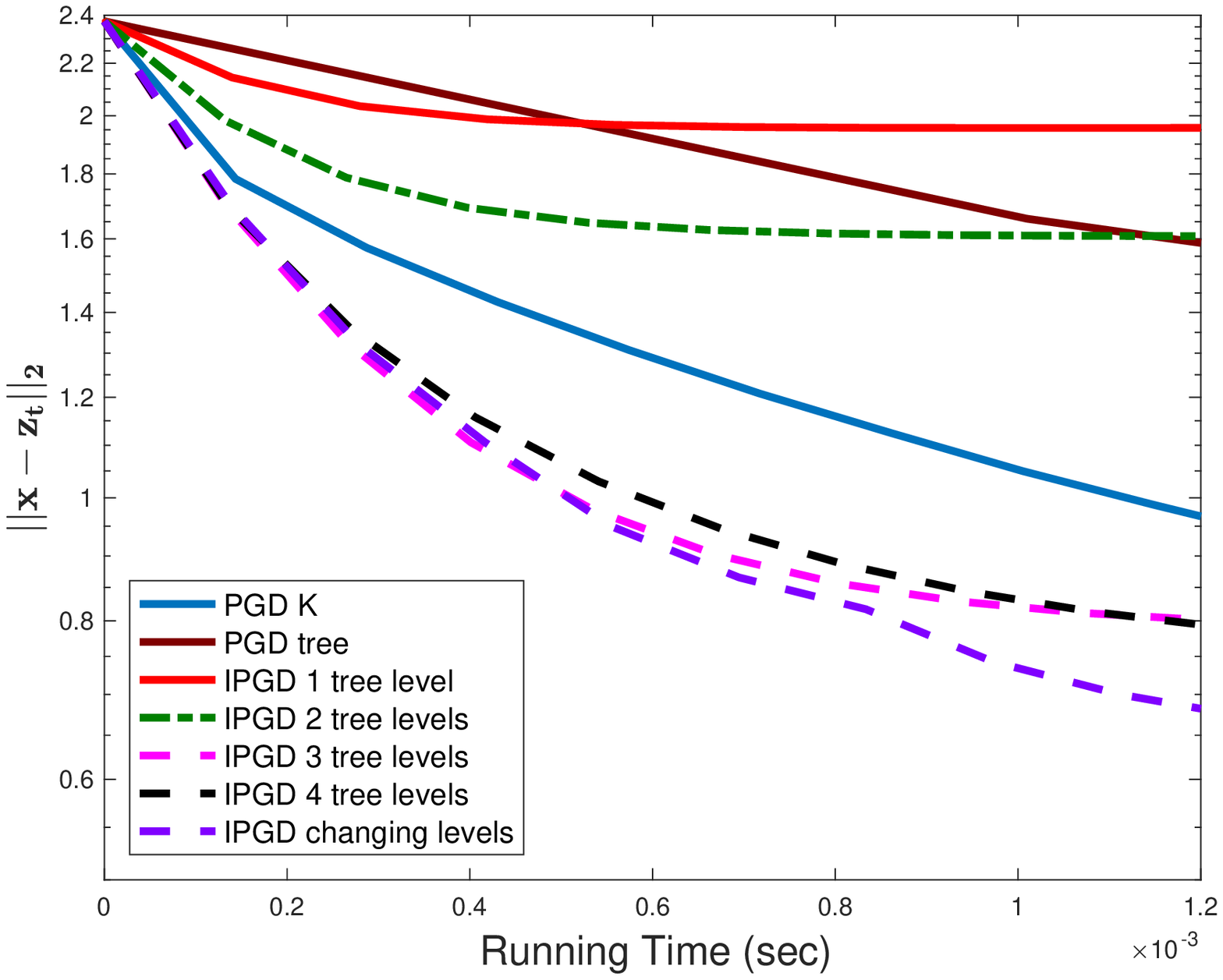}
\label{fig:tree_error_run_time_zoom}}
}
}
\end{center}
\caption{Reconstruction error as a function of the iterations (left) and the running time (right) for recovering a sparse vector with tree structure. 
Since we initialize all algorithms with the zero vector, the error at iteration/time zero is $\norm{\vect{x}}$. Zoomed version of the first $10$ iterations and first $1ms$ appears in the bottom row.
This figure demonstrates the convergence rate of PGD with projections onto the sparse set and  sparse tree set compared to IPGD with $p$ that projects onto a certain number of levels of the tree and IPGD with changing $p$ that projects onto an increasing number of levels as the iterations proceed. Note that while PGD with a projection onto a tree structure converges faster than IPGD as a function of the number of iterations (left figure), it converges slower than IPGD if we take into account the actual run time of each iteration, as shown in the right figure, due to the higher complexity of the PGD projections.  
}
\label{fig:tree}
\end{figure*}

The error term in theorems~\ref{thm:IPGD_error_t} and \ref{thm:IPGD_error_K} at iteration $t$ is comprised of two components. The first goes to zero as $t$ increases while the second increases with iterations and is on the order of $\epsilon$. The fewer iterations we perform the larger $\epsilon$ we may allow. An alternative perspective is that the larger the reconstruction error we can tolerate, the larger $\epsilon$ may be and thus we require fewer iterations. 
Therefore, the projection $p$ introduces a tradeoff. On the one hand, it leads to an increase in the reconstruction error. On the other hand, it simplifies the projected set, which leads to faster convergence (to a solution with larger error).

The works in \cite{Giryes14Greedy, Giryes15GreedySignal, Hegde14Approximation} use a similar concept of near-optimal projection (compared to \cite{Baraniuk10Model} that assumes only exact projections). The main difference between these contributions and ours is that these papers focus on specific models, while we present a general framework that is not specific to a certain low-dimensional prior. 
In addition, in these papers the projection is performed to make it possible to recover a vector from a certain low-dimensional set, while in this work the main purpose of our inexact projections is to accelerate the convergence within a limited number of iterations. For a larger number of iterations these projections may not lead to a good reconstruction error.

\section{Examples} 
\label{sec:examples} 
 
This section presents examples of IPGD with an operator $p$ that accelerates the convergence of PGD for a given set $\Kappa$.
 
\subsection{Sparse recovery with tree structure}
\label{sec:tree} 

To demonstrate our theory we consider a variant of the $k$-sparse set with tree structure in \eqref{eq:Kappa_hat_tree} that has smaller weights in the lower nodes of the tree. 
We generate a $k$-sparse vector $\vect{x} \in \RR{127}$ with $k = 13$ and a sparsity pattern that obeys a tree structure. Moreover, we generate the non-zero entries in $\vect{x}$ independently from a Gaussian distribution with  zero mean and variance $\sigma^2=1$ if they are at the first two levels of the tree %, $\sigma^2 = 0.5^2$ if they are at the third and fourth levels, 
 and $\sigma^2 = 0.2^2$ for the rest.

The best way to recover $\vect{x}$ is by using a projection onto the set $\hat\Kappa$ in \eqref{eq:Kappa_hat_tree}, which is the strategy proposed in the context of model-based compressed sensing \cite{Baraniuk10Model}. Yet, this projection requires some additional computations at each iteration \cite{Baraniuk10Model}. Our technique suggests to approximate it by a linear projection onto the first levels of the tree (a simple operation) followed by a projection onto $\Kappa = \{\vect{z} : \norm{\vect{z}}_0 \le k \}$.

The more levels we add in the projection $p$, the smaller the approximation error $\epsilon$ turns out to be. More specifically, it is easy to show that $\epsilon$ in \eqref{eq:P_K_pv_x_px_epsilon} is bounded by two times the energy of the entries eliminated from $\vect{x}$ divided by the total energy of $\vect{x}$, i.e., by $2\frac{\norm{p(\vect{x}) - \vect{x}}}{\norm{\vect{x}}}$. Clearly, the more layers we add the smaller $\epsilon$ becomes. Yet, assuming that all nodes in each layer are selected with equal probability, the probability of selecting a node at layer $l$ is equal to $\prod_{i=1}^l 0.5^{i-1}$, where we take into account the fact that a node can be selected only if all its forefathers have been chosen.  Thus, the upper layers have more significant impact on the values of $\epsilon$.

On the other hand, the convergence rate $\rho_p(\Kappa)$ for a projection with $l$ layers is equivalent to the convergence rate for the set of vectors of size $2^{l}$ (denoted by $\Kappa_l$). Thus, we get that $\rho_p(\Kappa) = \rho(\Kappa_l)$, which is dependent on the Gaussian mean width $\omega_{\Kappa_l}$ that scales as $\max(k l, k\log(2^l/k))$. Clearly, when we take all the layers $l = \log(d)$ and we have $\omega_{\Kappa_l} = \omega_\Kappa = O(k\log(d/k))$.

Figure~\ref{fig:tree_error} presents the signal reconstruction error ($\norm{\vect{x} - \vect{z}_t}_2$) as a function of the number of iterations for PGD with the sets $\Kappa$ (IHT \cite{Blumensath09Iterative}) and $\hat\Kappa$ (model-based IHT \cite{Baraniuk10Model})\footnote{For demonstration purposes we plot only the cases where model-based IHT converges to zero.} and for the proposed IPGD with $p$ that projects onto a different number of levels (1-5) of the tree. All algorithms use step size $\mu =\frac{1}{(\sqrt{d}+\sqrt{m})^2}$.
It is interesting to note that if $p$ projects only onto the first layer, then the algorithm does not converge as the resulting approximation error $\epsilon$ is too large. 
However, starting from the second layer, we get a faster convergence at the first iterations with $p$ that projects onto a smaller set, which yields a smaller $\rho$. As the number of iterations increases, the more accurate projections  achieve a lower reconstruction error, where the plateau attained is proportional to the approximation error of $p$ as predicted by our theory.
% Clearly, as PGD does not introduce an error in its projection, it achieves the smallest recovery error but not at the first iterations. This demonstrates the tradeoff between the approximation error we may allow and the  convergence speed. 

This tradeoff can be used to further accelerate the convergence by changing the projection in IPGD over the iterations. Thus, in the first iterations we enjoy the fast convergence of the coarser projections and in the later ones we use more accurate projections that allow achieving a lower plateau. The last line in Fig.~\ref{fig:tree} demonstrates this strategy, where at the first iteration $p$ is set to be a projection onto the first two levels,  and then every four iterations another tree level is added to the projection until it becomes a projection onto all the tree levels (in this case IPGD coincides with PGD). Note that IPGD converges faster than PGD also when the projection in it becomes onto all the tree levels. This can be explained by the fact that typically convergence of non-linear optimization techniques depends on the initialization point \cite{Bertsekas99Nonlinear}.

While here we arbitrarily chose to add another level every fixed number of iterations, in general, a control set can be used for setting the number of iterations to be performed in each training level. We demonstrate this strategy in Section~\ref{sec:LIPGD}.

Since PGD with $\hat{\Kappa}$ does not introduce an error in its projection and projects onto a precise set, it achieves the smallest recovery error throughout all iterations. Yet, as its projection is computationally demanding, it converges slower than IPGD if we take into account the run time of each iteration, as can been seen in Fig.~\ref{fig:tree_error_run_time}. This clearly demonstrates the advantage of using simple projections with IPGD compared to accurate but more complex projections with PGD.

\begin{figure}[t]
\begin{center}
{
{\subfigure[House image.]{\includegraphics[width=0.23\textwidth]{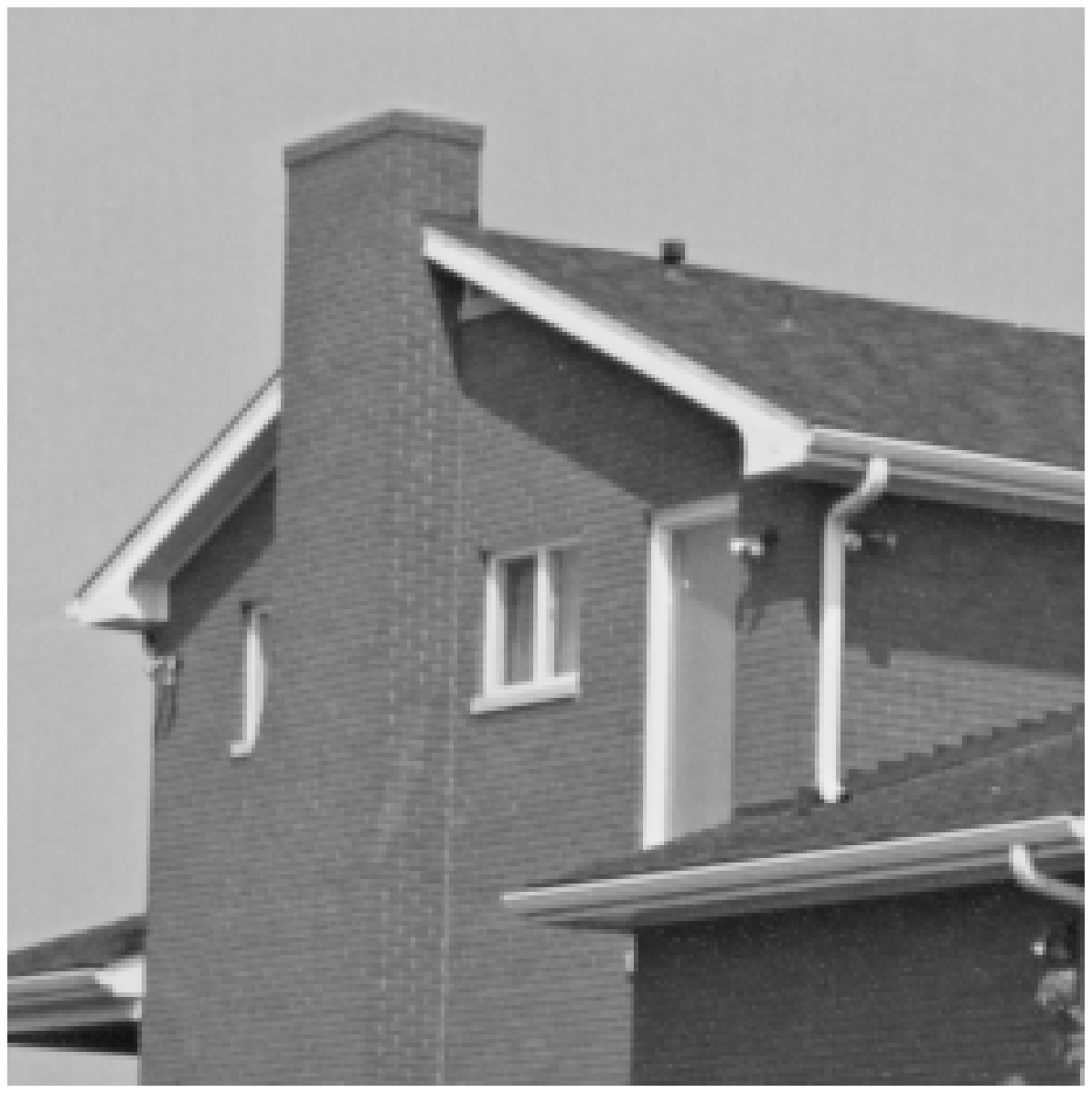}}
\label{fig:house_image}
\hspace{0.1in}
\subfigure[Patch from house image.]{\includegraphics[width=0.23\textwidth]{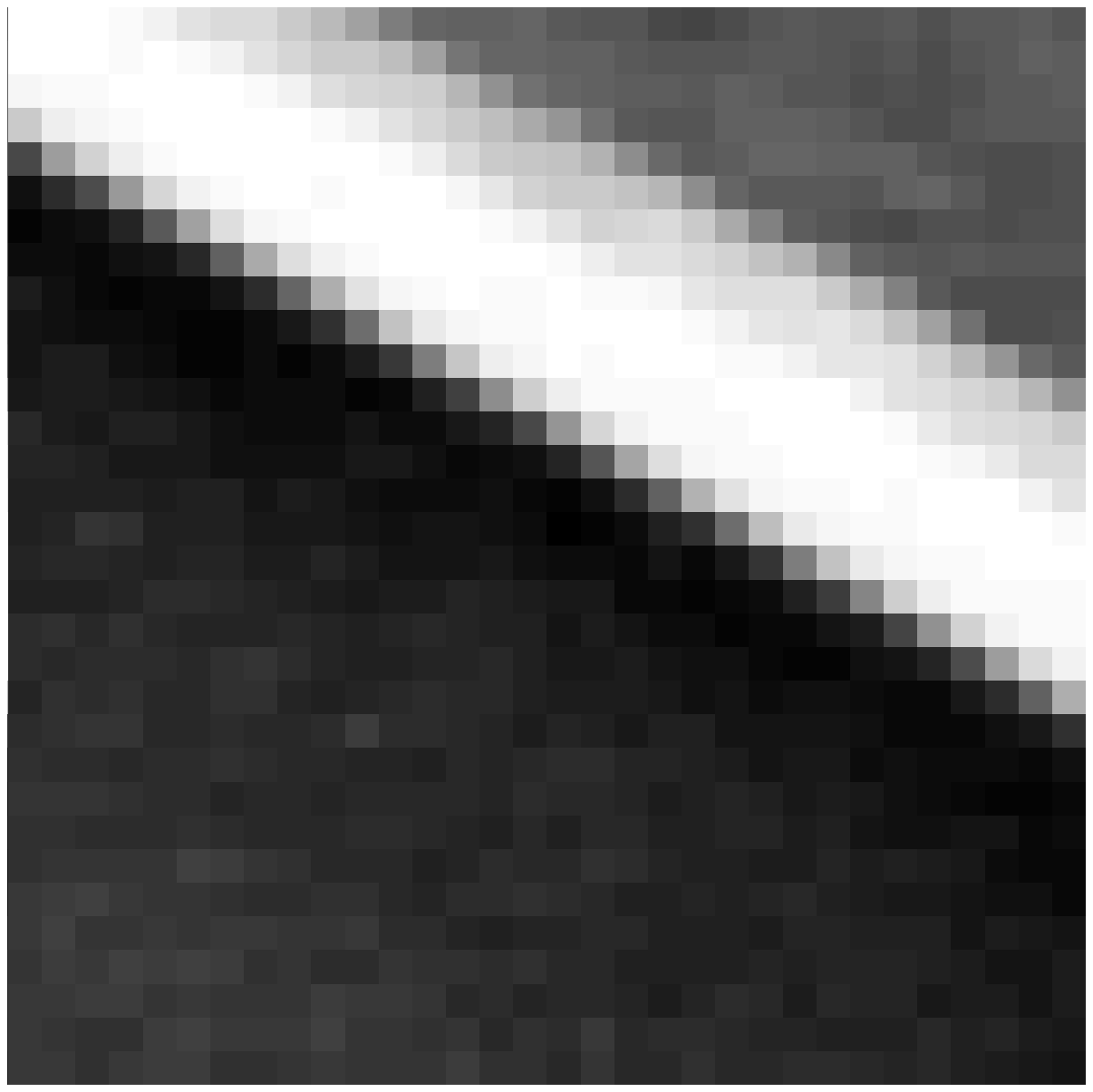}}
\label{fig:house_patch} \\  \vspace{0.2in}
\subfigure[Patch representation magnitudes with DCT.]{\includegraphics[width=0.23\textwidth]{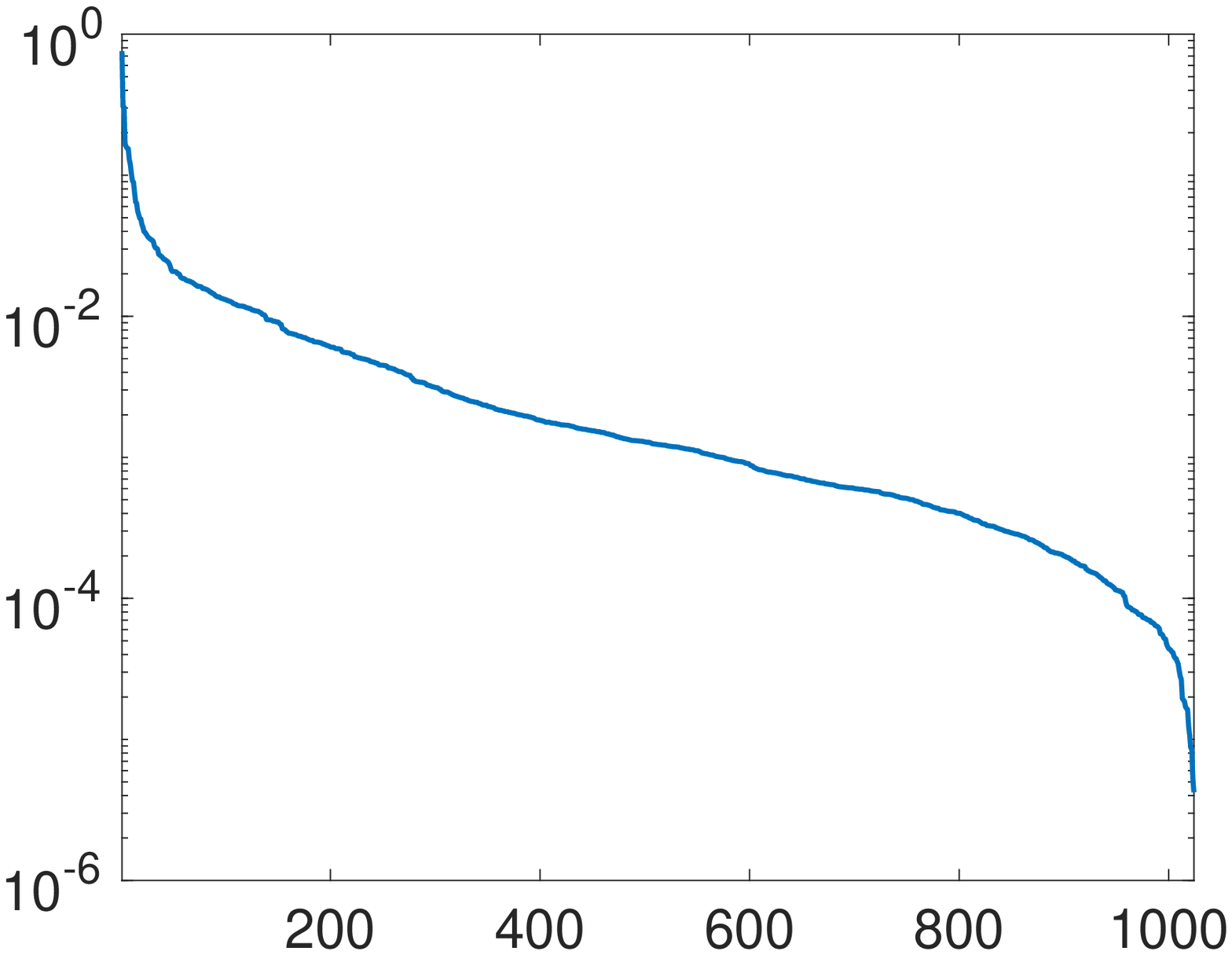}}
\label{fig:house_dct}
\hspace{0.05in}
\subfigure[Patch representation magnitudes with Haar.]{
\includegraphics[width=0.23\textwidth]{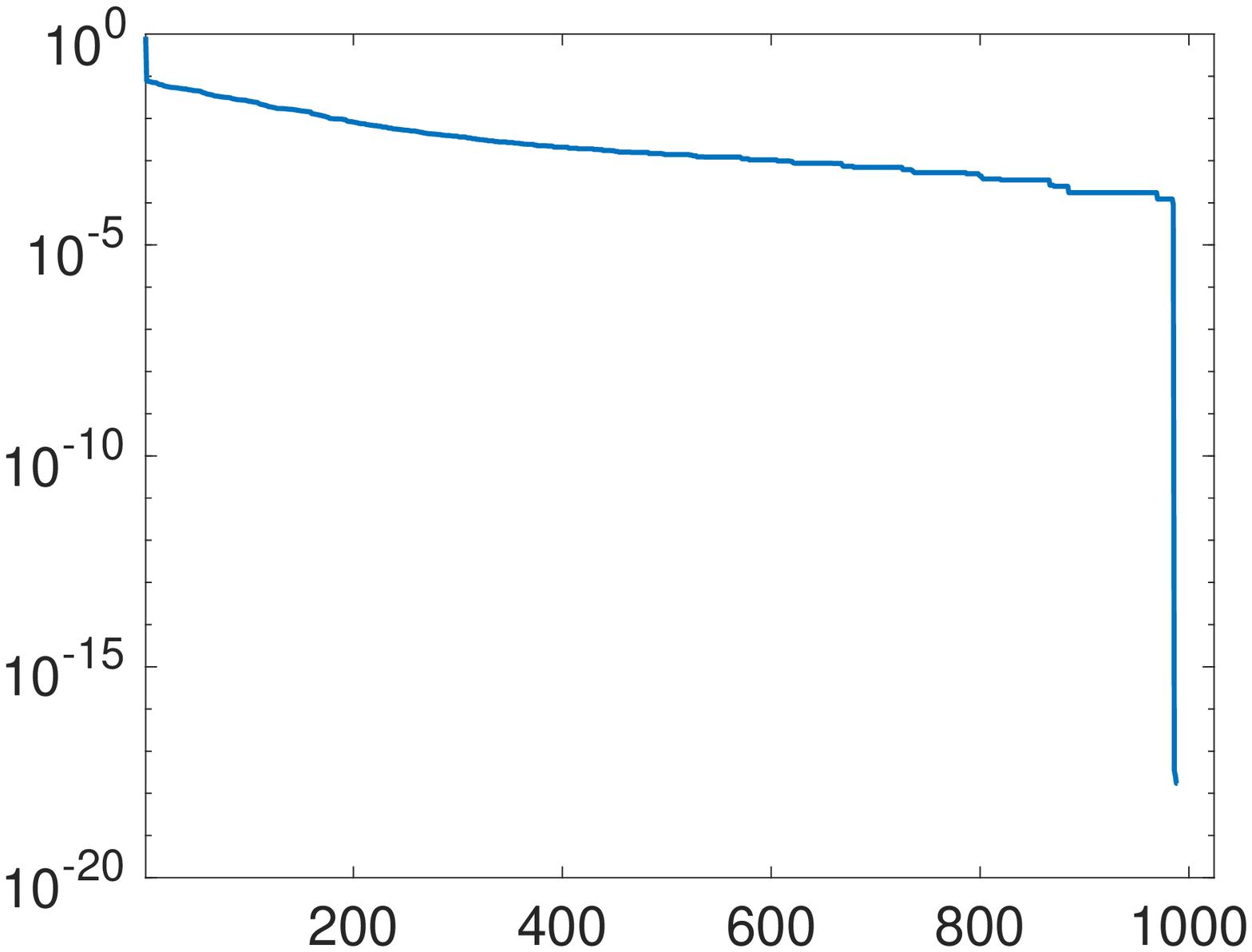}}
\label{fig:house_haar}
}
}
\end{center}
\caption{House Image (top left) and a random patch selected from it (top right) with the sorted magnitude (in log scale) of the representation of this patch in the DCT (bottom left) and Haar (bottom right) bases.}
\label{fig:house}
\end{figure}

\begin{figure}[t]
\begin{center}
{
{\includegraphics[width=0.48\textwidth]{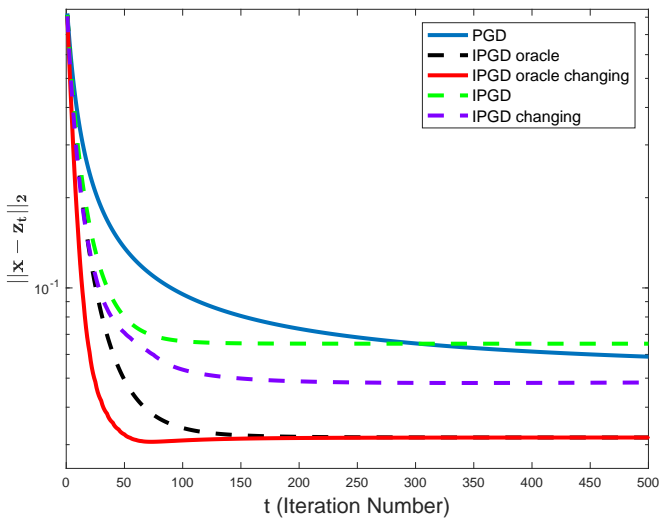}
\label{fig:side_inf_error}}
}
\end{center}
\caption{Reconstruction error as a function of the iterations for sparse recovery with side information. This demonstrates the convergence rate of (i) PGD with a projection onto the $\ell_1$ ball compared to (ii) IPGD with oracle side information on the columns of the representation of $\vect{x}$ in the Haar basis; (iii) IPGD with oracle side information that projects onto an increasing number of columns from the Haar basis ordered according to their significance in representing $\vect{x}$; (iv) IPGD with a projection onto the first $512$ columns of that Haar basis; and (v) IPGD with a changing $p$ that projects onto an increasing number of columns from the Haar basis.}
\label{fig:side_inf}
\end{figure}

\subsection{Sparse recovery with side information}
\label{sec:sparse_side}

Another possible strategy to improve reconstruction that  relates to our framework is using side information about the recovered signal, e.g., from estimates of similar signals. This approach was applied  to improve the quality of MRI \cite{Fessler92Regularized, Hero99Minimax, Weizman15Compressed} and CT scans \cite{Chen08Prior}, and also in the general context of sparse recovery  \cite{Friedlander12Recovering, Khajehnejad09Weighted, Vaswani10Modified, Weizman15Compressed}.

We  demonstrate this approach, in combination with our proposed framework, for the recovery of a sparse vector under the discrete cosine transform (DCT), given information of its representation under the Haar transform. Our sampling matrix is $\matr{M} = \matr{A}\matr{D}^*$, where $\matr{A} \in \RR{700 \times 1024}$ is a random matrix with i.i.d. normally distributed entries, $\matr{D}^*$ is the (unitary) DCT transform (that is applied on the signal before multiplying the DCT coefficients with the random matrix $\matr{A}$) and $\matr{D}$ is the DCT dictionary.  We use random patches of size $32\times 32$, normalized to have unit $\ell_2$ norm, from the standard {\em house} image. Note that such a patch is not exactly sparse either in the Haar or the DCT domains. See Fig.~\ref{fig:house} for an example of one patch of the house image.
Without considering the side information of the Haar transform, one may recover $\vect{x}$ (a representation of a patch in the DCT basis) by using PGD with the set $\Kappa = \{\vect{z} : \norm{\vect{z}}_1 \le \norm{\vect{x}}_1 \}$. Given $\vect{x}$ we may recover the patch $\matr{D}\vect{x}$.

Assume that someone gives us oracle side information on the set of Haar columns corresponding to the largest coefficients that contain $95\%$ of the energy in a patch $\matr{D}\vect{x}$.  While there are many ways to incorporate the side information in the recovery, we show here how IPGD can be used for this purpose.
Denoting by  $\matr{P}_{\vect{x}, 95\%}^{oracle}$ the linear projection onto this set of columns, one may apply IPGD with $p = \matr{D}\matr{P}_{\vect{x}, 95\%}^{oracle}\matr{D}^*$ and $\Kappa$. As $\norm{\matr{D}\vect{x}} = \norm{\vect{x}}$ (since $\matr{D}$ is unitary), we have that $\epsilon = 0.05$ in \eqref{eq:Kp_eps_ineq}.  Figure~\ref{fig:side_inf} compares between PGD with $\Kappa$ and IPGD with $p$ and $\Kappa$. We average over $100$ different randomly selected sensing matrices and patches.

The number of columns in Haar that contain $95\%$ of the energy is roughly $d/2$. Thus, the Gaussian mean width $\omega_p(C)$ in this case is roughly the width of the tangent cone of the $\ell_1$ norm at a $k$-sparse vector in the space of dimension $d/2$, which is smaller than $\omega(C)_{\ell_1}(\vect{x})$. Thus, $\rho_p(C)$ is smaller than $\rho(C)$.
Clearly, for less energy preserved in $\vect{x}$ (i.e., bigger $\epsilon$) we need less columns from Haar, which implies a smaller Gaussian mean width and a faster convergence rate $\rho_p(\Kappa)$. We have here a tradeoff between the approximation error we may allow $\epsilon$ and the convergence rate $\rho_p(\Kappa)$ that improves as $\epsilon$ increases.

Since projections onto smaller sets lead to faster convergence we suggest as in the previous example to apply PGD with an oracle projection that uses less columns from the Haar basis at the first iteration (i.e., has larger $\epsilon$) and then adds columns gradually throughout the iterations. The third (red) line in Fig.~\ref{fig:side_inf} demonstrates this option, where the first iterations use  a projection onto the columns that contain $50\%$ of the energy of $\vect{x}$ and then every $5$ iterations the next $50$ columns correspoding to the coefficients with the largest energy are added. We continue until the columns span $95\%$ of the energy of the signal. Thus, IPGD with changing projections converges faster than IPGD with a constant $p = \matr{D}\matr{P}_{\vect{x}, 95\%}^{oracle}\matr{D}^*$ but reaches the same plateau.  

Typically, oracle information on the coefficients of $\vect{x}$ in the Haar basis is not accessible. Even though, it is still possible to use common statistics of the data to accelerate convergence. For example, in our case it is known that most of the energy of the signal is concentrated in the low-resolution Haar filters. Therefore, we propose to use IPGD with a projection $p$ that projects onto the first $512$ columns of the Haar basis. As before, it is possible to accelerate convergence by projecting first on a smaller number of columns and then increasing the number as the iterations proceed (in this case we add columns till IPGD coincides with PGD). These two options are presented in the fourth and fifth line of Fig.~\ref{fig:side_inf}, respectively. Both of these options provide faster convergence, where IPGD with a fixed projection $p$ incurs a higher error as it uses less accurate projections in the last iterations compared to PGD and IPGD with changing projections. The plateau of the latter is the same one of the regular PGD (which is not attained in the graph due to its early stop) but is achieved with a much smaller number of iterations.

\begin{table}
\footnotesize
\begin{center}
\begin{tabular}{|c|c|c|c|c|}
\hline
Image & Bicubic & OMP & IHT & LIPGD \\
\hline
baboon & 23.2 & 23.5 & 23.4 & \bf 23.6\\
bridge & 24.4 & 25.0 & 24.8 & \bf 25.1 \\
coastguard & 26.6 & 27.1 & 26.9 & \bf 27.2 \\ 
comic & 23.1 & 24.0 & 23.8 & \bf 24.2  \\
face & 32.8 & 33.5 & 33.2 & \bf 33.6 \\
flowers & 27.2 & 28.4 & 28.1 & \bf 28.7 \\
foreman & 31.2 & 33.2 & 32.3 & \bf 33.5 \\
lenna & 31.7 & 33.0 & 32.6 & \bf 33.2 \\
man & 27.0 & 27.9 & 27.7 & \bf 28.1 \\
monarch & 29.4 & 31.1 & 30.9 & \bf 31.6 \\ 
pepper & 32.4 & 34.0 & 33.6 & \bf 34.4 \\
ppt3 & 23.7 & 25.2 & 24.6 & \bf 25.5 \\
zebra & 26.6 & 28.5 & 28.0 & \bf 28.9 \\
\hline
\end{tabular}
\end{center}
\caption{PSNR of super-resolution by bicubic interpolation and a pair of dictionaries with various sparse coding methods.}
\label{tbl:sup_res}
\end{table}

\section{Learning the Projection -- Learned IPGD (LIPGD)}
\label{sec:LIPGD}

\begin{figure}[t]
\begin{center}
{{\includegraphics[width=0.48\textwidth]{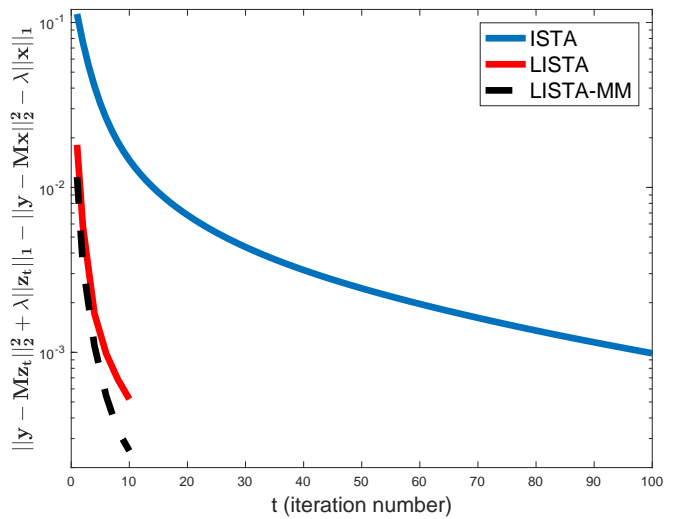}}
}
\end{center}
\caption{The $\ell_1$ loss as a function of the iterations of ISTA, LISTA and LISTA-MM applied on patches from the {\em house} image. This demonstrates the faster convergence of the proposed LISTA-MM compared to LISTA and the fast convergence of LISTA compared to ISTA.}
\label{fig:LISTA}
\end{figure}

In many scenarios, we may not know what type of simple operator $p$ causes $\Proj_\Kappa(p(\cdot))$ to approximate $\hat{\Kappa}$ in the best possible way. Therefore, a useful strategy is to learn $p$ for a given dataset. Assuming a linear $p$, we may rewrite \eqref{eq:IPGD} as
\begin{eqnarray}
\label{eq:IPGD_net}
\vect{z}_{t+1} = \Proj_{\Kappa}\left( p\left(\mu\matr{M}^*\vect{y}\right) + 
p\left(\left(\matr{I}  - \mu \matr{M}^*\matr{M}\right)\vect{z}_t \right)\right).
\end{eqnarray}
Instead of learning $p$ directly, we may learn two matrices $\matr{A}$ and $\matr{U}$, where the first replaces $p\mu\matr{M}^*$ and the second $p\left(\matr{I}  - \mu \matr{M}^*\matr{M}\right)$. 
This results in the iterations
\begin{eqnarray}
\label{eq:Learned_IPGD}
\vect{z}_{t+1} = \Proj_{\Kappa}\left( \matr{A}\vect{y} + 
\matr{U}\vect{z}_t \right),
\end{eqnarray}
which is very similar to those of LISTA in \eqref{eq:LISTA}. The only difference between \eqref{eq:Learned_IPGD} and LISTA is the non-linear part, which is an orthogonal projection in the first and a proximal mapping in the second. 
 
We apply this method to replace the sparse coding step in the super-resolution algorithm proposed in \cite{Zeyde12single}, where a pair of low and high resolution dictionaries is used to reconstruct the patches of the high-resolution image from the low-resolution one. In the code provided by the authors of \cite{Zeyde12single}, orthogonal matching pursuit (OMP) \cite{MallatZhang93} with sparsity $3$ is used. The complexity of this strategy corresponds to IHT with $3$ iterations. The target sparsity we use with IHT is higher ($k=40$) as it was observed to provide better reconstruction results. Note that in IHT, unlike OMP, the number of iterations may be different than the sparsity level. For optimal hyperparameter selection (such as choosing the target sparsity level), we use the training set used for the training of the dictionary in \cite{Zeyde12single}, which contains $91$ images.
 
Since IHT does not converge with only $3$ iterations, we apply LIPGD to accelerate convergence. We use the same dictionary dimension as in \cite{Zeyde12single} ($30 \times 1000$ and $81\times1000$ for the low and high resolution dictionaries, respectively), and train an LIPGD network to infer the sparse code of the image patches in the low-resolution dictionary. Training of the weights is performed by stochastic gradient descent with batch-size $1000$ and Nesterov momentum \cite{Nesterov83Method} for adaptively setting the learning rate. We train the network using only the first $85$ images in the training set, keeping the last $6$ as a validation set. We reduce the training rate by a factor of $2$ if the validation error stops decreasing. The initial learning rate is set to $0.001$ and the Nesterov parameter to $0.9$. 
We use the sparse representations of the training data calculated by IHT or LIPGD to generate the high-resolution dictionary as in \cite{Zeyde12single}.

Table~\ref{tbl:sup_res} summarizes the reconstruction results of regular bicubic interpolation, the OMP-based super-resolution technique of  \cite{Zeyde12single} (with $3$ iterations) and its version with IHT and LIPGD (replacing OMP). It can be seen clearly that IHT leads to inferior results compared to OMP since it does not converge in $3$ iterations. LIPGD improves over both IHT and OMP as the training of the network allows it to provide good sparse approximation with only $3$ iterations. This demonstrates the efficiency of the proposed LIPGD technique, which has the same computational complexity of both OMP and IHT.

\section{Learning the Projection -- LISTA Mixture Model}
\label{sec:cs_DL} 
 
Though the theory in this paper applies directly only to \eqref{eq:Learned_IPGD} (with some constraints on $\matr{A}$ and $\matr{U}$ that stem from the constraints on $p$), the fast convergence of LISTA may be explained by the resemblance of the two methods. The success of LISTA may be interpreted as learning to approximate the set $\hat{\Kappa}$ in an indirect way by learning the linear operators $\matr{A}$ and $\matr{U}$.
In other words, it can be viewed as a method for learning a linear operator that together with the proximal mapping $\shrink_{f, \lambda}$ approximates a more accurate proximal mapping of a true unknown function $\hat{f}$ that leads to much faster convergence.

With this understanding, we argue that using multiple inexact projections may lead to faster convergence as each can approximate in a more accurate way different parts of the set $\hat{\Kappa}$.  In order to show this, we propose a LISTA mixture model (MM), similar to the Gaussian mixture model proposed in \cite{Yu12Solving}, in which we train several LISTA networks, one for each part of the dataset. Then, once we get a new vector, we apply all the networks on it in parallel (and therefore with negligible impact on the latency, which is very important in many applications) and chose the one that attains the smallest value in the objective of the minimization problem \eqref{eq:min_unconst}. 

We test this strategy on the {\em house} image by extracting from it patches of size $5 \times 5$, adding random Gaussian noise to each of them with variance $\sigma^2= 25$ and then removing the DC and normalizing each. We take $7/9$ of the patches for training and $1/9$ for validation and testing. We train LISTA to minimize directly the objective \eqref{eq:min_unconst} as in \cite{Sprechmann15Learning} and stop the optimization after the error of the validation set increases.
For the LISTA-MM we use $6$ LISTA networks such that we train the first one on the whole data. We then remove $1/6$ of the data whose objective value in \eqref{eq:min_unconst} is the closest to the one ISTA attains after  $1000$ iterations. We use this LISTA network as the initialization of the next one that is trained on the rest of the data. We repeat this process by removing in the same way the part of the data with the smallest relative error and then train the next network. After training $6$ networks we cluster the data points by selecting for each patch the network that leads to the smallest objective error for it in \eqref{eq:min_unconst} and fine tune each network for its corresponding group of patches. We repeat this process $5$ times. The objective error of \eqref{eq:min_unconst} as a function of the number of iterations/depth of the networks is presented in Fig.~\ref{fig:LISTA}. Indeed, it can be seen that partitioning the data, which leads to a better approximation, accelerates  convergence. 

Our proposed LISTA-MM strategy bears some resemblence to the recently proposed rapid and accurate image super resolution (RAISR) algorithm \cite{Romano17RAISR}. In this method, different filters are trained for different types of patches in natural images. This leads to improved quality in the attained up-scaled images with only minor overhead in the computational cost, leading to a very efficient super-resolution technique.

\section{Conclusion}
\label{sec:conc}

In this work we suggested an approach to trade-off between approximation error and convergence speed. This is accomplished by approximating complicated projections by  inexact ones that are computationally efficient. We provided theory for the convergence of an iterative algorithm that uses such an approximate projection and showed that at the cost of an error in the projection one may achieve faster convergence in the first iterations. The larger the error the smaller the number of iterations that enjoy fast convergence. This suggests that if we have a budget for only a small number of iterations (with a given complexity), then it may be worthwhile to use inexact projections which can result in a worse solution in the long term but make better use of the given computational constraints. Moreover, we showed that even when we can afford a larger number of iterations, it may be worthwhile to use  inexact projections in the first iterations and then change to more accurate ones at latter stages. %which can lead to an improved reconstruction.

Our theory offers an explanation to the recent success of neural networks for 
approximating the solution of certain minimization problems. These networks achieve similar accuracy to iterative techniques developed for such problems (e.g., ISTA for $\ell_1$ minimization) but with much smaller computational cost. We demonstrate the usage of this method for the problem of image super-resolution. In addition, our analysis provides a technique for estimating the solution of these minimization problems by using multiple networks but with fewer layers in each of them.

\appendices

\section{Proof of Theorem~\ref{thm:IPGD_error_t}}
\label{sec:IPGD_error_t_proof}

The proof of Theorem~\ref{thm:IPGD_error_t} relies on the following lemma.

\begin{lem}
\label{lem:Pc_IMM_diff_bound}
Under the same conditions of Theorem~\ref{thm:IPGD_error_t}
\begin{eqnarray}
&& \hspace{-0.35in}\norm{\Proj_{C}\left(p\left( \matr{I} - \mu \matr{M}^*\matr{M}\right)(\vect{z}_t - \vect{x})\right)} \\ \nonumber && \hspace{0.8in} \le \rho_p(C)\norm{\vect{z}_t - \vect{x}} + \epsilon(1+\rho_p(C)) \norm{\vect{x}}.
\end{eqnarray}
\end{lem}

\noindent {\it Proof:}
Since $\vect{z}_t = \Proj_\Kappa(p \vect{v})$ for a certain vector $\vect{v}$, we have
\begin{eqnarray}
\label{eq:IPGD_lem_conv_proof_1_pre}
&& \hspace{-0.28in} \norm{\Proj_{C}\left(p\left( \matr{I} - \mu \matr{M}^*\matr{M}\right)(\vect{z}_t - \vect{x})\right)} 
\\ && \nonumber \hspace{-0.24in}
\le \sup_{\substack{\vect{v} \in \Real^d \text{ s.t. } \\   \norm{\Proj_{\Kappa}(p\vect{v})- \vect{x}} 
\le \norm{\vect{z}_t - \vect{x}}}} \norm{\Proj_{C}\left(p\left( \matr{I} - \mu \matr{M}^*\matr{M}\right)\left(\Proj_{\Kappa}(p\vect{v})- \vect{x} \right)\right)}
\\ && \nonumber \hspace{-0.24in}
= \sup_{\substack{\vect{v} \in \Real^d \text{ s.t. } \\   \norm{\Proj_D(p\vect{v}- \vect{x})} 
\le \norm{\vect{z}_t - \vect{x}}}} \norm{\Proj_{C}\left(p\left( \matr{I} - \mu \matr{M}^*\matr{M}\right)\Proj_D(p\vect{v} - \vect{x})\right)},
\end{eqnarray}
where the last equality follows from Lemma~\ref{lem:Kappa_D_rel}.
Using the triangle inequality with \eqref{eq:IPGD_lem_conv_proof_1_pre} leads to
\begin{eqnarray}
\label{eq:IPGD_lem_conv_proof_1}
&& \hspace{-0.32in} \norm{\Proj_{C}\left(p\left( \matr{I} - \mu \matr{M}^*\matr{M}\right)(\vect{z}_t - \vect{x})\right)} 
\\ && \nonumber \hspace{-0.28in} \le 
\sup_{\substack{\vect{v} \in \Real^d \text{ s.t. } \\   \norm{\Proj_D(p\vect{v}- \vect{x})} 
\le \norm{\vect{z}_t - \vect{x}}}} \hspace{-0.3in} \norm{\Proj_{C}\left(p\left( \matr{I} - \mu \matr{M}^*\matr{M}\right)\Proj_D(p\vect{v} - p\vect{x})\right)} 
\\ && \nonumber \hspace{-0.28in}
+ \norm{\Proj_{C}\left(p\left( \matr{I} - \mu \matr{M}^*\matr{M}\right)\frac{\Proj_D(\vect{x} - p\vect{x})}{\norm{\Proj_D(\vect{x} - p\vect{x})}}\right)}
 \norm{\Proj_D(p\vect{x} - \vect{x})}.
\end{eqnarray}

We turn now to bound the first and second terms in the right-hand-side (rhs) of \eqref{eq:IPGD_lem_conv_proof_1}. For the second term, note that since  $\Proj_D(\vect{x} - p\vect{x}) \in D$, we have
\begin{eqnarray}
\label{eq:IPGD_lem_conv_subproof_rho}
&& \hspace{-0.6in}\norm{\Proj_{C}\left(p\left( \matr{I} - \mu \matr{M}^*\matr{M}\right)\frac{\Proj_D(\vect{x} - p\vect{x})}{\norm{\Proj_D(\vect{x} - p\vect{x})}}\right)} \\
&& \nonumber \le \sup_{\vect{v} \in D \cap \mathbb{B}^d} \norm{\Proj_{C}\left(p\left( \matr{I} - \mu \matr{M}^*\matr{M}\right)\vect{v}\right)}  \\
&& \nonumber \le \sup_{\vect{v}, \vect{u} \in C \cap \mathbb{B}^d} \norm{\vect{u}^*p\left( \matr{I} - \mu \matr{M}^*\matr{M}\right)\vect{v}} \le \rho(\Kappa),
\end{eqnarray}
where the second inequality follows from Lemma~\ref{lem:P_c_sup} and the fact that $D \subset C$. In the last inequality we replace $\vect{u}^*p = (p^*\vect{u})^*$ by $\tilde{\vect{u}}$ and take the supremum over it instead of over $\vect{u}$. In addition, $ \norm{\Proj_D(p\vect{x} - \vect{x})} \le \epsilon$ from \eqref{eq:Dp_eps_ineq}. 

For the first term in the rhs of \eqref{eq:IPGD_lem_conv_proof_1} we use the inverse triangle inequality $\norm{\Proj_D(p\vect{v}- p\vect{x})} - \norm{\Proj_D(p\vect{x}- \vect{x}) }  \le \norm{\Proj_D(p\vect{v}- \vect{x})}$ and \eqref{eq:Dp_eps_ineq}. Combining the results leads to
\begin{eqnarray}
\label{eq:IPGD_step_error_proof_2_convex}
&& \hspace{-0.3in} \norm{\Proj_{C}\left(p\left( \matr{I} - \mu \matr{M}^*\matr{M}\right)(\vect{z}_t - \vect{x})\right)}  
 \\ && \nonumber \hspace{-0.3in}
\le \hspace{-0.14in} \sup_{\substack{\vect{v} \in \Real^d \text{ s.t. } \\   \norm{\Proj_D(p(\vect{v}- \vect{x}))} 
\le \norm{\vect{z}_t - \vect{x}}+ \epsilon\norm{\vect{x}}}} \hspace{-0.35in} \norm{\Proj_{C}\left(p\left( \matr{I} - \mu \matr{M}^*\matr{M}\right)\Proj_D(p(\vect{v} - \vect{x}))\right)} 
\\ && \nonumber \hspace{2.4in}
+ \epsilon \rho(\Kappa)   \norm{\vect{x}}
 \\ && \nonumber \hspace{-0.3in}
=\hspace{-0.14in} \sup_{\substack{\vect{u}\in C\cap \mathbb{B}^d,
\vect{v} \in \Real^d  \\  \text{ s.t. }\norm{\Proj_D(p\vect{v})} 
\le 1}} \hspace{-0.15in}
\norm{\vect{u}^*p\left( \matr{I} - \mu \matr{M}^*\matr{M}\right)\Proj_D(p\vect{v})}\left(\norm{\vect{z}_t - \vect{x}}+ \epsilon\norm{\vect{x}}\right)
\\ && \nonumber \hspace{2.4in}
+ \epsilon \rho(\Kappa)   \norm{\vect{x}}
\\ && \nonumber \hspace{-0.3in}
=\hspace{-0.14in} \sup_{\vect{u} \in C\cap \mathbb{B}^d, p\vect{v} \in D\cap \mathbb{B}^d} \hspace{-0.15in}
\norm{(p\vect{u})^*\left( \matr{I} - \mu \matr{M}^*\matr{M}\right)p\vect{v}}\left(\norm{\vect{z}_t - \vect{x}}+ \epsilon\norm{\vect{x}}\right)
\\ && \nonumber \hspace{2.4in}
+ \epsilon \rho(\Kappa)   \norm{\vect{x}}
\\ && \nonumber \hspace{-0.3in}
\le \hspace{-0.14in} \sup_{\vect{v},\vect{u} \in C\cap \mathbb{B}^d} 
\norm{(p\vect{u})^*\left( \matr{I} - \mu \matr{M}^*\matr{M}\right)p\vect{v}}\left(\norm{\vect{z}_t - \vect{x}}+ \epsilon\norm{\vect{x}}\right)
\\ && \nonumber \hspace{2.4in}
+ \epsilon \rho(\Kappa)   \norm{\vect{x}},
\end{eqnarray}
where the first equality follows from Lemma~\ref{lem:P_c_sup} and the second from the definition of a projection onto a set. The last inequality follows from the fact that $\{ p\vect{v} \in D \cap \mathbb{B}^d \} \subset D \cap \mathbb{B}^d \subset C \cap \mathbb{B}^d$. Reordering the terms and using the definition of $\rho_p(C)$ leads to the desired result.
 \hfill $\Box$

 We turn now to the proof of Theorem~\ref{thm:IPGD_error_t}

{\it Proof:} The IPGD error at iteration $t+1$ is,
\begin{eqnarray}
\label{eq:IPGD_error_def}
 \norm{\vect{z}_{t+1} - \vect{x}} = \norm{\Proj_{\Kappa}\left(p\left(\vect{z}_t + \mu \matr{M}^*(\vect{y} - \matr{M}\vect{z}_t)\right)\right) - \vect{x}}.
\end{eqnarray}
Using Lemma~\ref{lem:Kappa_D_rel} and the fact that $\vect{y} = \matr{M}\vect{x}$ we have
\begin{eqnarray}
\label{eq:IPGD_step_error_proof_ab}
&& \hspace{-0.33in} \norm{\vect{z}_{t+1} - \vect{x}}   \hspace{-0.025in}
= \hspace{-0.025in} \norm{\Proj_{D}\left(p\left( \matr{I} - \mu \matr{M}^*\matr{M}\right)(\vect{z}_t - \vect{x}) - \vect{x} + p\vect{x}\right) }
\\ \nonumber &&  \overset{(a)}{\le}  \norm{\Proj_{D}\left(p\left( \matr{I} - \mu \matr{M}^*\matr{M}\right)(\vect{z}_t - \vect{x})\right) } 
+ \norm{\Proj_{D}\left(\vect{x} - p\vect{x} \right) }
\\ \nonumber &&  \overset{(b)}{\le}   \norm{\Proj_{C}\left(p\left( \matr{I} - \mu \matr{M}^*\matr{M}\right)(\vect{z}_t - \vect{x})\right)} 
+ \epsilon\norm{\vect{x}},
\end{eqnarray}
where $(a)$ follows from the convexity of $D$  and the triangle inequality; and $(b)$ from \eqref{eq:Dp_eps_ineq} and Lemma~\ref{lem:C_D_rel}. Using Lemma~\ref{lem:Pc_IMM_diff_bound} with  \eqref{eq:IPGD_step_error_proof_ab} leads to
\begin{eqnarray}
\label{eq:IPGD_step_error_proof_iter}
\norm{\vect{z}_{t+1} - \vect{x}}  \le \rho_p(\Kappa)\norm{\vect{z}_t - \vect{x}} + \epsilon(2+\rho_p(\Kappa)) \norm{\vect{x}}.
\end{eqnarray}
Applying the inequality in \eqref{eq:IPGD_step_error_proof_iter} recursively provides the desired result.  
 \hfill $\Box$ 

\section{Proof of Theorem~\ref{thm:IPGD_error_K}}
\label{sec:IPGD_error_K_proof}

The proof Theorem~\ref{thm:IPGD_error_K} relies on the following lemma.

\begin{lem}
\label{lem:PK_IMM_diff_bound}
Under the same conditions of Theorem~\ref{thm:IPGD_error_K},
\begin{eqnarray}
&& \hspace{-0.25in}\norm{\Proj_{\Kappa-\Kappa}\left(p\left( \matr{I} - \mu \matr{M}^*\matr{M}\right)(\vect{z}_t - \vect{x})\right)} \\ \nonumber && \hspace{0.45in} \le \rho_p(\Kappa)\norm{\vect{z}_t - \vect{x}} + \epsilon(2\rho(\Kappa)+\rho_p(\Kappa)) \norm{\vect{x}}.
\end{eqnarray}
\end{lem}

{\it Proof:} 
Using Lemma~\ref{lem:Kappa_D_rel} and the fact that $\vect{z}_t = \Proj_\Kappa (p\vect{v})$ for a certain vector $\vect{v} \in \RR{d}$ and then the triangle inequality, leads to
\begin{eqnarray}
\label{eq:IPGD_lemma_error_proof_1_nonconvex}
&& \hspace{-0.3in} \norm{\Proj_{\Kappa-\Kappa}\left(p\left( \matr{I} - \mu \matr{M}^*\matr{M}\right)(\vect{z}_t - \vect{x})\right)}  
\\ && \nonumber \hspace{-0.3in}
\le\sup_{\substack{\vect{v} \in \Real^d \text{ s.t. } \\   \norm{\Proj_D(p\vect{v}- \vect{x})} 
\le \norm{\vect{z}_t - \vect{x}}}} \hspace{-0.3in} \norm{\Proj_{\Kappa-\Kappa}\left(p\left( \matr{I} - \mu \matr{M}^*\matr{M}\right)\Proj_D(p\vect{v} - \vect{x})\right)} 
\\ && \nonumber \hspace{-0.3in}
\le\sup_{\substack{\vect{v} \in \Real^d \text{ s.t. } \\   \norm{\Proj_D(p\vect{v}- \vect{x})} 
\le \norm{\vect{z}_t - \vect{x}}}} \hspace{-0.3in} \norm{\Proj_{\Kappa-\Kappa}\left(p\left( \matr{I} - \mu \matr{M}^*\matr{M}\right)\Proj_D(p\vect{v} - p\vect{x})\right)} +
\\ && \nonumber \hspace{-0.3in}
\norm{\Proj_{\Kappa-\Kappa}\left(p\left( \matr{I} - \mu \matr{M}^*\matr{M}\right)\frac{\Proj_D(p\vect{v} - \vect{x}) - \Proj_D(p\vect{v} - p\vect{x})}{\norm{\Proj_D(p\vect{v} - \vect{x}) - \Proj_D(p\vect{v} - p\vect{x})}}\right)}
\\  &&  \nonumber \hspace{1.25in}
 \cdot \norm{\Proj_D(p\vect{v} - \vect{x}) - \Proj_D(p\vect{v} - p\vect{x})}.
\end{eqnarray}
Using \eqref{eq:P_D_pv_x_pv_px_epsilon} and the same steps of the proof of Theorem~\ref{thm:PGD_error_K}, we may bound the second term in the rhs of \eqref{eq:IPGD_lemma_error_proof_1_nonconvex} by $2\epsilon \rho(\Kappa)   \norm{\vect{x}}$. This leads to 
\begin{eqnarray}
\label{eq:IPGD_lemma_error_proof_2_nonconvex_pre}
&& \hspace{-0.3in} \norm{\Proj_{\Kappa-\Kappa}\left(p\left( \matr{I} - \mu \matr{M}^*\matr{M}\right)(\vect{z}_t - \vect{x})\right)}  
\\ && \nonumber \hspace{-0.3in}
\le\sup_{\substack{\vect{v} \in \Real^d \text{ s.t. } \\   \norm{\Proj_D(p\vect{v}- \vect{x})} 
\le \norm{\vect{z}_t - \vect{x}}}} \hspace{-0.3in} \norm{\Proj_{\Kappa-\Kappa}\left(p\left( \matr{I} - \mu \matr{M}^*\matr{M}\right)\Proj_D(p\vect{v} - p\vect{x})\right)}
 \\ && \nonumber \hspace{1.465in}
+ 2\epsilon \rho(\Kappa)   \norm{\vect{x}}.
\end{eqnarray}
From the inverse triangle inequality together with \eqref{eq:P_D_pv_x_pv_px_epsilon}, we have that $\norm{\Proj_D(p\vect{v}- \vect{x})} \ge \norm{\Proj_D(p\vect{v}- p\vect{x})} - \epsilon \norm{\vect{x}}$. Thus, 
\begin{eqnarray}
\label{eq:IPGD_lemma_error_proof_2_nonconvex}
&& \hspace{-0.3in} \norm{\Proj_{\Kappa-\Kappa}\left(p\left( \matr{I} - \mu \matr{M}^*\matr{M}\right)(\vect{z}_t - \vect{x})\right)}  
 \\ && \nonumber \hspace{-0.2in}
\le \hspace{-0.23in} \sup_{\substack{\vect{v} \in \Real^d \text{ s.t. } \\   \norm{\Proj_D(p(\vect{v}- \vect{x}))} 
\le \norm{\vect{z}_t - \vect{x}}+ \epsilon\norm{\vect{x}}}} \hspace{-0.4in} \norm{\Proj_{\Kappa-\Kappa}\left(p\left( \matr{I} - \mu \matr{M}^*\matr{M}\right)\Proj_D(p(\vect{v} - \vect{x}))\right)} 
\\ && \nonumber \hspace{1.465in}
+ 2\epsilon \rho(\Kappa)   \norm{\vect{x}}
\\ && \nonumber \hspace{-0.2in}
\le 
%\hspace{-0.13in} \sup_{\vect{v},\vect{u} \in (\Kappa-\Kappa)\cap \mathbb{B}^d} \hspace{-0.1in} \norm{(p\vect{u})^*\left( \matr{I} - \mu \matr{M}^*\matr{M}\right)p\vect{v}}
\rho_p(\Kappa)\left(\norm{\vect{z}_t - \vect{x}}+ \epsilon\norm{\vect{x}}\right) 
+ 2\epsilon \rho(\Kappa)   \norm{\vect{x}},
\end{eqnarray}
where the last inequality follows from the same line of argument  used for deriving \eqref{eq:IPGD_step_error_proof_2_convex} in Lemma~\ref{lem:Pc_IMM_diff_bound} (with $\Kappa-\Kappa$ instead of $C$). 
 \hfill $\Box$ 

We now turn to the proof of Theorem~\ref{thm:IPGD_error_K}.

{\it Proof:} Denoting $\tilde{\vect{v}} = \left( \matr{I} - \mu \matr{M}^*\matr{M}\right)(\vect{z}_t - \vect{x})$, the IPGD error at iteration $t+1$ obeys
\begin{eqnarray}
\label{eq:IPGD_step_error_proof_dc}
&& \hspace{-0.5in} \norm{\vect{z}_{t+1} - \vect{x}}  
=  \norm{\Proj_{D}\left(p\tilde{\vect{v}} - \vect{x} +p \vect{x}\right) }
\\ \nonumber && \hspace{-0.2in} 
\overset{(c)}{\le}  \norm{\Proj_{D}\left(p\tilde{\vect{v}} \right) }
+  \norm{\Proj_{D}\left(p\tilde{\vect{v}}- \vect{x} +p \vect{x}\right) - \Proj_{D}\left(p\tilde{\vect{v}} \right) }
\\ \nonumber && \hspace{-0.2in} 
\overset{(d)}{\le}  \kappa_\Kappa\norm{\Proj_{K-K}\left(p\tilde{\vect{v}} \right) }
+  \epsilon\norm{\vect{x} },
\end{eqnarray}
where $(c)$ follows from the triangle inequality; and $(d)$ from \eqref{eq:P_D_pv_x_px_epsilon} and Lemma~\ref{lem:C_D_rel}.
Using Lemma~\ref{lem:PK_IMM_diff_bound} with \eqref{eq:IPGD_step_error_proof_dc}, we get 
\begin{eqnarray}
\label{eq:IPGD_step_error_proof_non_convex}
&& \hspace{-0.25in} \norm{\vect{z}_{t+1} - \vect{x}}  
\\ \nonumber && \hspace{-0.1in} \le  \rho_p(\Kappa)\kappa_\Kappa\norm{\vect{z}_t - \vect{x}} + \epsilon(2\rho(\Kappa)\kappa_\Kappa+\rho_p(\Kappa)\kappa_\Kappa+1) \norm{\vect{x}}.
\end{eqnarray}
Applying \eqref{eq:IPGD_step_error_proof_non_convex} recursively leads to the desired result.  
 \hfill $\Box$ 

\section*{Acknowledgments}

RG is partially supported by GIF grant no. I-2432-406.10/2016 and ERC-StG grant no. 757497 (SPADE). YE is partially supported by the ERC grant no. 646804-ERC-COG-BNYQ. AB is partially supported by ERC-StG RAPID. GS is partially supported by ONR, NSF, NGA, and ARO.
We thank Dr. Pablo Sprechmann for early work and insights into this line of research, and Prof. Ron Kimmel and Prof. Gilles Blanchard for insightful comments.

{\small
\bibliographystyle{plain}
\bibliography{references}
}

\end{document}